\documentclass[reqno,11pt]{amsart}
\usepackage{amsmath, amssymb, amsthm, amsfonts} 
\usepackage[english]{babel}
\usepackage{bbm}
\usepackage{graphicx}
\usepackage{url}
\usepackage{array} 
\usepackage{endnotes}
\usepackage{epstopdf}
\usepackage{subcaption}
\usepackage[a4paper,bindingoffset=0.5cm,left=2cm,right=2cm,top=2.5cm,bottom=2cm,footskip=.8cm]{geometry}

\usepackage{pgfplotstable}
\usepgfplotslibrary{groupplots}
\usepackage{caption}
\usepgfplotslibrary{statistics}

\newcommand{\nac}[1]{\textcolor{red}{ #1}}

\newcommand{\plotPair}[5]{%
  \addplot+[smooth, thick, mark=*, color=#4, mark options={color=#4, fill=#4}]
    table[col sep=comma, x=i, y=#2]{#1};
  \addplot+[smooth, thick, dashed, mark=x, color=#4, mark options={color=#4}]
    table[col sep=comma, x=i, y=#3]{#1};
}

\newcommand{\plotAbsPair}[5]{%
  \addplot+[smooth, thick, mark=*, mark size=#5, color=#4,
            mark options={color=#4, fill=#4}]
    table[col sep=comma, x=n, y=#2]{#1};
  \addplot+[smooth, thick, dashed, mark=x, mark size=#5, color=#4,
            mark options={color=#4}]
    table[col sep=comma, x=n, y=#3]{#1};
}

\usepackage{tabularx}
\usepackage{hyperref}
\hypersetup{
    colorlinks=true,
    linkcolor=blue,
    filecolor=blue,      
    urlcolor=blue,
}

\usepackage{tikz}
\usetikzlibrary{arrows, automata,positioning,calc,shapes,decorations.pathreplacing,decorations.markings,shapes.misc,petri,topaths}
  \usepackage{pgfplots}
  \pgfplotsset{compat=newest}
  \usetikzlibrary{plotmarks}
  \usepackage{grffile}
\newlength\figureheight
  \newlength\figurewidth
\pgfplotsset{%
    tick label style={font=\scriptsize},
    label style={font=\footnotesize},
    legend style={font=\footnotesize},
         every axis plot/.append style={very thick}
}
\usetikzlibrary{math,shapes.geometric}

\usepackage{rotating}
\usepackage{amsbsy,enumerate}
\usepackage{graphicx}
\usepackage{comment}
\usepackage{mathrsfs} 
\usepackage{diagbox}
\usepackage{floatrow}
\usepackage{longtable}
\usepackage{booktabs}
\usepackage[ruled,vlined]{algorithm2e}
\renewcommand{\SetArgSty}[1]{%
  \renewcommand{\ArgSty}[1]{
    \textnormal{\csname#1\endcsname{##1}}\unskip
  }%
}

\definecolor{bcolA}{HTML}{440154} 
\definecolor{bcolB}{HTML}{30678d} 
\definecolor{bcolC}{HTML}{35b778} 
\definecolor{bcolD}{HTML}{fde724}
\definecolor{bcolA1}{HTML}{440154} 
\definecolor{bcolB1}{HTML}{20908C} 
\definecolor{bcolC1}{HTML}{FDE724} 

\newcommand{\bs}{\boldsymbol}

\newcommand{\vb}{\vspace{3.2mm}}

\allowdisplaybreaks

\newtheorem{lemma}{Lemma}

\newtheorem{theorem}{Theorem}
\newtheorem{remark}{Remark}

\newtheorem{definition}{Definition}

\usepackage{abstract}
\begin{document}


\title[Window design in delivery schedules]{A framework for window design in delivery schedules}

\author{R. Bekker, B. Bharti, N. Levering, M. Mandjes}

\maketitle

\begin{abstract} This paper develops a structured framework for the design and dynamic updating of service time windows in delivery and appointment-based systems. We consider a single-server setting with stochastic service and travel times, where customers are promised a time window in which the provider will arrive. The first part of the paper introduces a static window construction method based on a probabilistic threshold criterion, using an analytical approximation of residual travel and service time distributions. Building on this, we develop a dynamic update mechanism that monitors residual system uncertainty, where time windows are revised during execution only when the remaining time until the window's start falls below a predefined threshold. This threshold-based approach enables communication-efficient scheduling while substantially improving delivery accuracy.
Numerical experiments demonstrate significant performance gains of the dynamic approach in both stylized and real-world settings.

\vb



\vb

\noindent
{\sc Affiliations.} 
BB is with the Korteweg-de Vries Institute for Mathematics, University of Amsterdam, Science Park 904, 1098 XH Amsterdam, The Netherlands.

\noindent RB  is with the
Department of Mathematics,
Vrije Universiteit Amsterdam,
De Boelelaan 1111,
1081 HV Amsterdam,
the Netherlands. 
RB is the corresponding author: r.bekker@vu.nl

\noindent During the writing
of this paper, NL was with the SARA team at LAAS, 7 Avenue du Colonel Roche, 31077 Toulouse, France; the RMESS team at IRIT – ENSEEIHT,
2 rue C. Camichel,                    
31071 Toulouse, France.
She is currently at the Netherlands Defense Academy, Den Helder, The Netherlands.

\noindent MM is with the Mathematical Institute, Leiden University, P.O. Box 9512,
2300 RA Leiden,
The Netherlands. He is also affiliated with Korteweg-de Vries Institute for Mathematics, University of Amsterdam, Amsterdam, The Netherlands; E{\sc urandom}, Eindhoven University of Technology, Eindhoven, The Netherlands; Amsterdam Business School, Faculty of Economics and Business, University of Amsterdam, Amsterdam, The Netherlands. 

\vb

\noindent 
{\sc Funding.}
NL acknowledges the support of 
Labex CIMI [Grant ANR-11-LABX-0040].
The research of BB and MM was supported by the European Union’s Horizon 2020 research and innovation programme under the Marie Skłodowska-Curie grant agreement no.\ 945045, and by the NWO Gravitation project NETWORKS under grant no.\ 024.002.003. \includegraphics[height=1em]{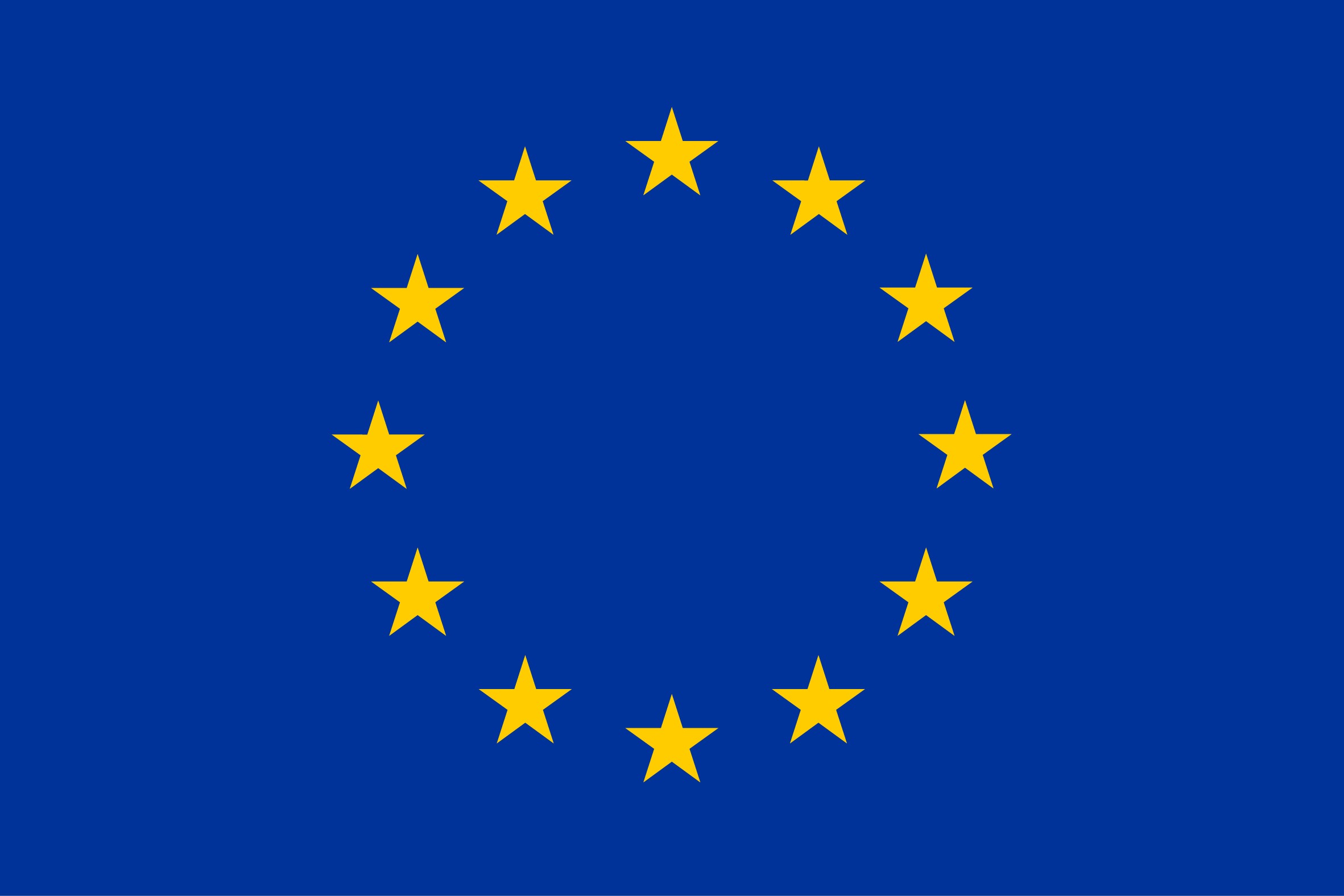}

\vb

\noindent
{\sc Acknowledgments.} 
The authors would like to thank Roshan Mahes for inspiring discussions. 

\end{abstract}

\maketitle

\section{Introduction}


The global parcel delivery market has seen rapid growth in recent years, driven by the surge in e-commerce and rising demand for fast and flexible delivery. Parcel volume rose from 64 billion in 2016 to 185 billion in 2023, with projections reaching 217 billion by 2025; in the U.S. alone, shipments exceeded 21 billion in 2022 (\url{https://www.contimod.com/package-statistics/}). This growth has triggered major investments in infrastructure (such as parcel lockers) and innovations such as real-time tracking, thus positioning the sector for continued expansion. Moreover, the growth has resulted in the rise of various companies that are specialized in parcel delivery, and that can be used by other companies to outsource the delivery component of their e-commerce. 

One way firms with a parcel delivery component aim to differentiate themselves is through the accuracy of their delivery forecasts. Customers are usually provided with time windows during which they are expected to be at a given location (typically their home) to receive their parcel. These time windows must be communicated well in advance, often the day before, so that customers can plan their day accordingly. When designing time windows, firms must account for various uncertainties. That is, a delivery driver visits a series of locations, each location corresponding to a customer. The total delivery time is influenced by two key sources of variability: the time required to complete each handoff (e.g., how long it takes for previous customers to answer the door) and the travel times between stops (which can be affected by factors such as traffic congestion and weather conditions). Both elements are inherently unpredictable, adding complexity to the task of reliable time window forecasting.

In recent years, the design of delivery time windows has taken on greater importance in e-commerce, largely due to the visibility and influence of online customer reviews. Since these reviews are easily shared and widely accessible, customer satisfaction has a direct impact on public perception—and, by extension, the future customer base—of a delivery firm. Customer satisfaction is heavily affected by two factors: (a) the reliability of the delivery (i.e., whether the parcel arrives within the specified window) and (b) the width of the delivery window (where wider windows are generally perceived more negatively). Importantly, there is an inherent trade-off between these two factors: achieving both high accuracy and narrow delivery windows can be difficult, given the uncertainty of the delivery operations. Wider time windows typically increase the likelihood of on-time delivery, but also require customers to remain available for longer periods.

In its simplest form, the window design problem can be described as follows: There are $n$ customers to be visited at known locations in a specified order. The travel times between locations are random, though their distribution is known; this time, in this paper referred to as `travel time', also includes the time spent handing over the parcel (e.g., parking, ringing the doorbell, having the customer sign to acknowledge receipt, etc.). 
We consider two scheduling scenarios: \textit{static} schedules, where time windows are set in advance and remain fixed, and \textit{dynamic} schedules, where time windows are updated in real-time based on the actual arrival times of earlier deliveries.
The primary objective of this paper is to present fast and intuitive methods for determining delivery windows, aiming to balance the reliability of the communicated time windows with their width.
By emphasizing speed and conceptual simplicity in window design, our approach is well-suited for practical deployment by service providers. Furthermore, the low computational cost makes our techniques suitable for integration with other systems, such as routing algorithms.

We focus on an objective function 
whose form is common in the appointment scheduling literature, as it combines three components of `disutility': it is a weighted sum of the expected times the deliverer arrives too early, the expected times they arrive too late, and a (not necessarily linear) function of the width of the time window. 
We show that, under mild assumptions, 
this objective allows for a per-customer decomposition, meaning it can be split into $n$ separate optimization problems. These assumptions include: (i)~if the deliverer arrives at the customer's location before the time window and the customer is not home, they do not wait until the window starts, (ii)~the service time does not depend on whether the parcel is actually handed over, and~(iii) the delivery company has not adopted a policy where, because of `fairness' reasons, the windows for individual customers must be equally wide.

The time window design problem has been analyzed in several papers; see the literature review in Section~\ref{sec:lit} below. 
In particular, in the static setting,
the objective function studied in
\cite{hosseini2025} represents a special
case of the one addressed in this paper.
Another recent paper closely aligned with our objectives is \cite{UGM}, where the goal is to minimize each customer's window width, subject to the constraint that the delivery occurs within the window with probability at least $1-\alpha$ (for some given `unreliability tolerance' $\alpha$). 
Notably, by incorporating penalties for both early and late arrivals, our model captures the practical reality that the probability of successfully delivering a parcel decreases as the deliverer deviates from the scheduled window, which is an effect that is not captured by the criterion used by \cite{UGM}.




{\it Contributions.} 
The first contribution of our work is a fast and intuitive method for designing delivery time windows along a given route, balancing both reliability and window width. For each customer on the route, this approach requires solving only a simple one-dimensional root-finding problem. 
The key equation involves the distribution of total travel time from the depot to the customer. We demonstrate that this arrival time can be effectively approximated by a suitably chosen normal random variable, enabling a highly efficient computational procedure. The accuracy of this approximation is thoroughly validated through numerical experiments. 
The presence of the arrival time distribution in the root-finding problem arises from our objective function, which minimizes the expected deviation between the actual delivery time and the communicated time window.

The second contribution examines the policy of offering all customers the {\it same} delivery window width, a practice often motivated by perceived fairness. While uniform window widths are commonly used by delivery companies as a means of competitive differentiation, imposing this constraint fundamentally alters the problem structure: the design problem no longer decouples, making it significantly more complex than the version without this restriction
(as addressed in our first contribution). Our experiments reveal that, despite its intention, the uniform policy actually does not lead to equitable outcomes: although all customers receive windows of equal width, those located earlier in the delivery route typically benefit from a higher likelihood of being served within the promised window. As an aside, we note that in the setup analyzed in \cite{UGM}, there is fairness in the sense that all customers are assigned the same unreliability tolerance $\alpha$, but there is unfairness in that some customers face substantially wider time windows than others.

The third contribution focuses on {\it dynamic schemes}, i.e., schemes in which customers receive window updates while the deliverer is in the process of delivering the $n$ parcels. With the advent of technologies that enable real-time notifications, it is now possible to provide customers with increasingly accurate information, allowing the communication of a narrower or shifted time window
during the delivery process. This means that if the deliverer falls behind schedule, they may choose to delay the windows of later customers, while if they are ahead of schedule, there is an incentive to bring the windows forward, thereby reducing the likelihood of arriving outside the communicated window. We demonstrate that a significant improvement can already be achieved by the dynamic
use of the static procedure, in which we carefully
select the moment at which a new window is communicated to the customer in order to leave
enough time between the communication time and
the time window, while simultaneously preventing a communication overload.
Our numerical experiments
show that such
a substantial improvement
exists for both synthetic 
homogeneous and heterogeneous scenarios, as well as for
settings that are 
based on a real-world delivery data set.

\smallskip

\textit{Outline.} This paper is organized as follows. Section~\ref{sec:lit} reviews the relevant literature. In Section~\ref{sec:framework}, we present our model and introduce the static and dynamic window assignment problems. The simplest case, the static window design problem without any uniformity requirement on the window widths, is addressed in Section~\ref{sec:wos}. This optimization problem is referred to as WOS (window-based optimal schedule) throughout the paper. In Section~\ref{sec:uwos}, we discuss the more computationally demanding case where all windows must be equally wide, which we term UWOS (window-based optimal uniform schedule). Finally, Section~\ref{sec:dyn} explores the potential gains that can be achieved by updating the schedule `on the fly', and Section~\ref{sec: conclusive_remarks} contains conclusive remarks.

\section{Literature review}\label{sec:lit}
In this section, we review relevant work on window design. Given the vast literature on delivery planning, our overview is not exhaustive; for additional references and further background, see, for example, \cite[\S2]{UGM} and \cite[\S2]{VRT}. Since our approach draws on concepts from the field of appointment scheduling, we also discuss its connection to that domain.

\medskip



{\it Time window assignment.} The works that are probably the most closely related to ours are Ulmer {\it et al.}~\cite{UGM}, Vareias {\it et al.}~\cite{VRT}, and Hosseini {\it et al.}~\cite{hosseini2025}. These papers consider what Jabali et al.~\cite{JAB} have termed a {\it self-imposed time window problem}, where the delivery company designs and communicates the time windows to the customers. This contrasts with the literature on attended home delivery, in which customers communicate
a time window or choose a time window from a presented set. There, the goal is to determine
a cost-minimizing vehicle routing schedule for this set of customer-selected time windows.
Recent papers that study variants of this problem include \cite{BKCE, KEC, PRLD}; 
for a comprehensive overview and an all-encompassing categorization of problem types, we refer to the survey on attended home delivery \cite{WAS}.
It should be noted that by allowing customers to choose their own time windows, one
limits the use of efficient routing procedures, thereby increasing the cost 
for delivery companies significantly.

In this paper, we adopt a setting — following, for example, \cite{UGM} — in which the routing and time window assignment decisions are {\it decoupled}: the route is assumed to be predetermined in an ‘outer loop’, and the focus is on assigning delivery time windows along this fixed route. This modeling choice is motivated in part by practical considerations — for instance, routing is often managed by third-party logistics providers — as noted in \cite{UGM}. Additionally, decoupling offers substantial computational advantages: the sequence of visits can be determined using a standard deterministic traveling salesman problem, after which the time window assignment is optimized for that route. While this approach may entail some loss in efficiency compared to fully integrated methods, the loss is often modest; see \cite{bekker2023} for a quantification of this effect in a related context. 

By contrast, a large body of literature addresses the {\it integrated} optimization of routing and time window assignment, resulting in more complex formulations in which the route and time windows are jointly optimized in the presence of uncertainty. These can be viewed as extensions of the traveling salesman problem under uncertainty \cite{LLM}, which does not capture the effects of uncertainties on customer satisfaction. Several studies \cite{dalmeijer2018, NPGAA, spliet2015} adopt an integrated approach assuming deterministic travel times and uncertain customer demands, with a capacity-constrained delivery vehicle. Jabali \textit{et al.}~\cite{JAB} address a similar problem but reverse the assumptions: demand is known while travel times are random, with each arc having a probability of disruption. However, their model restricts disruptions to at most one arc. Their solution method is two-staged — first determining the route, then scheduling along that fixed route.
This model is extended by Vareias \textit{et al.}~\cite{VRT}, who allow multiple arc disruptions, each with independently distributed continuous durations. They also propose a second model in which arc travel times are represented as independent discrete random variables, yielding a scenario-based formulation. This discrete model includes penalties for earliness, lateness, and window width, as well as a reliability constraint for on-time arrivals. The second-stage problem is an LP in the continuous case and a MIP in the discrete case.
Recent work has extended the integrated routing and time window design framework to settings with uncertain travel times. 
Similar to the above works,
Hosseini \textit{et al.}~\cite{hosseini2025} present 
a two-stage solution-method, where
their window design is based on a three-term criterion with a linear penalty component for the window width.
Hosseini \textit{et al.}~\cite{hosseini2025} also focus on
the setting in which the underlying
travel-time distributions are not known.
Hoogeboom \textit{et al.}~\cite{HADJ} address a similar challenge, thereby using robust optimization techniques.

Similar to our work, Ulmer {\it et al.}~\cite{UGM} purely considers the problem of time window assignment. However, their study focuses on assigning time windows to clients at the moment a request is made, rather than when all tour locations are known. This assignment is based on the system state~$s$ at the time of the request, which may include information about clients already on the route and a probabilistic model of future requests.
The objective is to minimize the expected width of the assigned time windows, where the expectation is taken over all possible states $s$. This is subject to a constraint on the probability that a client is served outside their assigned window --- again averaged over all states, weighted by their respective probabilities.
In contrast, our approach employs an objective function that explicitly accounts for early and late arrivals, as well as the width of the time window. 
A key advantage is that it captures not only whether the deliverer arrives too early or too late, but also by how much. Ulmer \textit{et al.}~\cite{UGM} illustrate their method using an example involving $10^4$ system states, with travel-time distributions estimated via Monte Carlo simulation for each. Our approach differs in relying on a simpler proxy, using primarily the mean and variance of travel-time distributions between successive clients. More broadly, a major strength of our method lies in its low computational overhead.

The works discussed above focus solely on the static time window assignment problem, where each customer is assigned a time window only once: either after all customer requests have been received or immediately upon request.
Celik \textit{et al.}~\cite{CSMW} explore a hybrid approach, in which the time window must be communicated between the request time and a pre-specified cut-off point, with customer satisfaction depending on the timing of this communication.
Similar
to our dynamic counterpart, 
Yu \textit{et al.}~\cite{YSBS} and Dalmeijer \textit{et al.}~\cite{DSW}
consider
a setting in which the communicated time windows can be updated during the 
delivery trip. In \cite{YSBS}, the authors propose re-optimizing initially assigned windows in response to new information; however, their model does not account for the impact of such updates on customer satisfaction.
This is taken into consideration in \cite{DSW},
who study a dynamic setting with two potential types of adjustment: time window extension
(time window start is fixed, width can be increased) and time window
postponement (width is fixed).

\medskip

{\it Appointment scheduling.} There are interesting connections between the window design problem addressed in this paper and the field of appointment scheduling. In the field of appointment scheduling, the objective function typically captures the aggregate `disutility' of all customers in the system; this is in essence also the approach we follow in the present paper. This objective function is then minimized over the vector of appointment times (rather than appointment windows). A common objective function might include the sum of the service provider's expected idle times (arising from completing a customer's service before the next appointment) and the customers' expected waiting times. Some examples of papers in this domain include \cite{hassin2008scheduling, kemper2014, PR, wang1997}, where it is noted that in all of them the order of service is assumed to be predetermined. 

The appointment problem in which the order is also optimized is intrinsically more challenging; see, for example, the overview paper \cite{ahmadi2017}. For an assessment of the intuitively appealing `smallest variance first' rule, see \cite{de2021}. In \cite{bekker2023, zhan2021}, integrated routing and appointment scheduling problems are analyzed. A specific feature of these studies is that if the deliverer arrives at the customer's location before the appointment time, they wait until this scheduled time. This makes these papers well-suited for applications such as home services, but less applicable to parcel delivery.

\section{Framework}\label{sec:framework}
The goal of this section is to define our model and introduce the optimization problems we aim to solve. The two sections that follow discuss the static and dynamic settings, respectively.


Let, for a given number $n\in {\mathbb N}$,  ${\bs B}\equiv (B_1,\ldots,B_n)$ be a non-negative random vector, where $B_i$ denotes the travel time of the $i$-th client ($i\in\{1,\ldots,n\}$). Here `travel time' is to be interpreted in a broad sense: e.g.\ in the delivery context $B_i$ could be the travel time from the location of client $i-1$ to that of client $i$, increased by the time it takes to hand over the parcel to client $i-1$ (where $B_1$ represents just the travel time from the depot to the first client).

The client is given a {\it window} during which she has to be at the delivery location. Let $(t_i,t_i+\Delta_i)$ be the window pertaining the $i$-th client. The deliverer does not wait for the client: in case the deliverer arrives at location $i$ before $t_i$ while the client is not present, the deliverer continues her route. This means that the partial sum
\[S_i:=\sum_{j=1}^i B_j\]
denotes the time the deliverer arrives at the location of client $i$, for $i=1,\ldots,n$. We denote by $F_i(t)$, for $t\in{\mathbb R}_+$, the cumulative distribution function ${\mathbb P}(S_i\leqslant t)$, where we impose the natural assumption that $F_i(0)=0$.

\begin{remark}\label{Remark: correlation}
{\em In most papers on window design that work with stochastic travel times, these are assumed independent, a notable exception being \cite{hassin2008scheduling}.
Observe that in our framework we do {\it not} impose the requirement that the components of the vector ${\bs B}$ be independent. This generality offers a substantial amount of flexibility. For instance, we could have that ${\bs B}= \bar {\bs B}\,X$, for a vector $\bar {\bs B}\in{\mathbb R}_+^n$ with independent entries, and $X\in {\mathbb R}_+$ encoding the impact of the weather or some other external factor. One could thus work with rules such as: when it rains, all travel times are 10\% longer; cf.\ also the framework adopted in \cite{UGM}.}
\hfill $\Diamond$
\end{remark}

\subsection{Description of the static problem}
In our most basic framework, the objective function to be minimized is, for $\omega \in (0,1)$,
\[G_n({\bs t},{\bs \Delta}):= \sum_{i=1}^n H_i({t}_i,{\Delta}_i),\quad
H_i(t,{\Delta}):=
\omega {\mathbb E}\left(S_i-t-\Delta\right)^+
+(1-\omega) {\mathbb E}\left(t-S_i\right)^+ +{\mathscr P}(\Delta),\]
to be minimized over ${\bs t}, {\bs\Delta}\in {\mathscr U}_n,$
where
\[{\mathscr U}_n:=\big\{{\bs t}, {\bs\Delta}\in {\mathbb R}_+^n:\:t_{i+1}\geqslant t_i\geqslant 0\:\mbox{\rm for all}\: i\in\{1,\ldots,n-1\}\big\}.\]
Here, $t_i$ is to be interpreted as the start time of the time window of client $i$ and $\Delta_i$ as the width of the time window.
The choice of this objective function is motivated as follows:
\begin{itemize}
\item[$\circ$]
The first term represents the cost corresponding to the scenario in which the deliverer is late, i.e. arrives at the client's location after $t_i+\Delta_i$, entailing that the client has to wait for the deliverer. The `loss' is proportional to the amount of time the deliverer is too late. 
\item[$\circ$]
The second term represents the cost corresponding to the scenario in which the deliverer is early, i.e. arrives at the client's location before $t_i.$ The `loss' being proportional to $t_i-S_i$ reflects that the more the deliverer is too early, the lower the chance of the client being present at the agreed location.
\item[$\circ$] The third term represents the penalty corresponding to the width of the window: evidently, ${\mathscr P}(\cdot)$ should be a positive non-decreasing function. 
\end{itemize}
This type of objective function is intensively worked with in the context of appointment scheduling: in that branch of the literature, the typical objective function consists of a weighted sum of the mean waiting times and mean idle times of the individual clients.

\begin{remark}\label{R1}
{\em
One could in principle work with more general objective functions, quantifying an asymmetry in risk aversion. A general form, with $\ell_i, \ell_i' \in \mathbb{R}_{+}$, is
\begin{align*}
    G_n({\bs t},{\bs \Delta}, {\bs \ell},{\bs \ell'}):= \omega\sum_{i=1}^n {\mathbb E}\left(\left(S_i-t_i-\Delta_i\right)^+\right)^{\ell_i}
+(1-\omega)\sum_{i=1}^n {\mathbb E}\left(\left(t_i-S_i\right)^+\right)^{\ell_i'} +\sum_{i=1}^n {\mathscr P}(\Delta_i).
\end{align*}
Observe that for larger values of $\ell_i$ and $\ell_i'$ the corresponding term is relatively flat, meaning that it has a relatively low impact on the value of the objective function. 
In the sequel, we exclusively consider the case of $\ell_i=\ell_i'=1$.}
\hfill $\Diamond$
\end{remark}

In the present subsection, we distinguish between two types of static optimization problems. In the first, the widths of the windows may vary. In the second, we impose the condition that all windows be of equal size. This latter formulation accommodates scenarios where maintaining consistent service intervals is required for operational or contractual reasons; cf.\ \cite[Section 2.1.1]{hosseini2025}.
We define
\[{\mathscr V}_n:=\big\{{\bs t}\in {\mathbb R}_+^n, {\Delta}\in {\mathbb R}_+:\:t_{i+1}\geqslant t_i\geqslant 0\:\mbox{\rm for all}\: i\in\{1,\ldots,n-1\}\big\}.\]

\begin{definition}[Window-based optimal schedule] For a given $n\in{\mathbb N}$, random vector ${\bs B}\in {\mathbb R}^n_+$, weight $\omega \in (0,1)$, and increasing penalty function ${\mathscr P}:{\mathbb R}_+\to {\mathbb R}_+$, define
\[{\rm WOS}_n({\bs B},\omega,{\mathscr P}):= \min_{{\bs t}, {\bs\Delta}\in {\mathscr U}_n}G_n({\bs t},{\bs \Delta}).\]
\end{definition}

\begin{definition}[Window-based optimal uniform schedule] For a given $n\in{\mathbb N}$, random vector ${\bs B}\in {\mathbb R}^n_+$, weight $\omega \in (0,1)$, and increasing penalty function ${\mathscr P}:{\mathbb R}_+\to {\mathbb R}_+$, define
\[{\rm UWOS}_n({\bs B},\omega,{\mathscr P}):= \min_{{\bs t},{\Delta}\in{\mathscr V}_n}G_n({\bs t},{\bs 1}\Delta).\]
\end{definition}

\begin{remark}\label{Remark: discrete time}
\normalfont
In practice, it may be preferable to communicate time windows with boundaries that are rounded to a multiple of 5 minutes, rather than using time values from a continuous range. However, following, for example, \cite{UGM,VRT}, we focus on the continuous setting and do not include the constraint that the resulting time windows have boundaries belonging to a pre-specified discrete set of values $\Gamma$. Note that this assumption is not restrictive, as for a relatively coarse 
$\Gamma$, the optimization procedure would simply involve computing the objective function value for all windows in 
$\Gamma$ and selecting the best option. Alternatively, for a fine 
$\Gamma$, the most natural procedure would be to execute a continuous-time procedure (i.e., the window assignment we intend to compute), after which proper rounding would yield the estimated optimal intervals in 
$\Gamma$.
\hfill $\Diamond$
\end{remark}



\subsection{Description of the dynamic problem}\label{subsec:dyn}
Up to this point, the time-window assignment has been performed only once, prior to the delivery tour -- a strategy that we have been referring to as {\it static} assignment. As outlined in the introduction, this paper also addresses the {\it dynamic} case, in which time windows are updated on the fly during the execution of the delivery route.

In practice, it is a desirable property not to provide the clients with excessively many update notifications. Supposing one has a cap on the number of updates, one could in principle use dynamic programming to compute the optimal window assignment at any point in time, using knowledge of the index of the current customer being `served' and the corresponding elapsed travel time. However, the dimensionality of the resulting optimization problem is prohibitively large, rendering such an approach impractical in operational settings.

Motivated by this limitation, we propose a more tractable alternative, referred to as DWOS (dynamic window-based optimal schedule): repeatedly applying the static solution obtained through WOS. As will be shown in Section \ref{sec:wos}, WOS admits a highly accurate approximation that can be computed in near real-time, making it well-suited for repeated use in a dynamic setting. The precise procedure for DWOS will be described in detail in Section \ref{sec:dyn}.

A key feature of the proposed approach is its limited communication requirement: each client is informed of a delivery time prior to the start of the delivery tour, and at most once more. The potential second notification occurs if the (dynamically updated) optimal start of the client's time window lies within $T$ time units of the current time, where $T$ is a parameter specified by the service provider (e.g., one or two hours). The rationale behind this limited notification policy is that frequent updates are typically perceived negatively by clients. An extension of this approach could involve offering clients real-time access to their current optimal delivery window via a mobile app, thereby maintaining transparency without over-communicating.


\section{Optimal static windows: the WOS case}\label{sec:wos}
In this section we present our analysis of the static WOS case. In Section \ref{subsec:foc} we present the underlying first-order conditions, and show that under mild conditions they have a unique solution. As it turns out, in case the travel times $B_i$, $i=1,\ldots,n$, stem from a normal distribution, a fairly explicit analysis can be performed; see Section \ref{subsec:WOSnorm}. The case with generally distributed travel times is addressed in Section \ref{subsec:WOSgen}.  

\subsection{First-order conditions} \label{subsec:foc}
We start by presenting the first-order conditions pertaining to the optimization problem 
${\rm WOS}_n({\bs B},\omega,\mathscr{P})$. Using standard differentiation rules, it is readily verified that these read
\begin{align}
    \frac{\partial}{\partial t_i} G_n({\bs t},{\bs \Delta}) &=  (1-\omega)F_i(t_i)- \omega(1-F_i(t_i + \Delta_i))=0 ,\label{FOC1}\\
    \frac{\partial}{\partial \Delta_i} G_n({\bs t},{\bs \Delta}) &= - \omega(1-F_i(t_i + \Delta_i)) + {\mathscr P}'(\Delta_i)=0 \label{FOC2}.
\end{align}
As mentioned above, the $2n$-dimensional optimization problem decouples into $n$ individual two-dimensional optimization problems that can be solved separately. This means that, in this static WOS setting, determining the optimal window assignment is of low computational complexity.

To gain some intuition about the nature of the solution of Equations \eqref{FOC1}--\eqref{FOC2}, we consider the specific case of a linear cost function:
${\mathscr P}(\Delta) = \alpha \Delta$. 
For any general underlying travel-time distribution and $\omega \in (0,1)$
we can obtain the corresponding result.
That is, combining the first two equations yields
\begin{align*}
    (1-\omega)F_i(t_i) = \alpha,
\end{align*}
from which we obtain that we should have $\alpha \leqslant 1 - \omega$ to have a stationary point. The second equation yields
\begin{align*}
    \omega(1-F_i(t_i + \Delta_i)) = \alpha,
\end{align*}
which only has a solution if $\alpha \leqslant \omega$. Therefore, only in the case
$\alpha \leqslant \min\{\omega, 1-\omega\}$ one has that the
optimum corresponds to a stationary point of the objective function (which implies that we should necessarily have $\alpha\leqslant \frac{1}{2}$).
Moreover, it is trivially seen that when
the inverse of the distribution function of
$S_i$ exists, one has
\begin{align*}
    t_i &= F_i^{-1}\left(\frac{\alpha}{1-\omega}\right),\quad
    \Delta_i = F_i^{-1}\left(1-\frac{\alpha}{\omega}\right) - F_i^{-1}\left(\frac{\alpha}{1-\omega}\right),
\end{align*}
in other words, the window starts at $F_i^{-1}(\frac{\alpha}{1-\omega})$ and ends at $F_i^{-1}(1-\frac{\alpha}{\omega})$.

The fact that we only get a non-degenerate result for $\alpha \leqslant \min\{\omega, 1-\omega\}$ is a consequence of the linear window cost: if ${\mathscr P}(\cdot)$ is proportional to a constant $\alpha$, then ${\mathscr P}'(\Delta)$ equals this $\alpha$, thus rendering solving the two first-order conditions simultaneously problematic. This complication does not play any role when picking, for instance, ${\mathscr P}(\Delta) = (\alpha/\beta) \Delta^\beta$, with $\alpha > 0$, $\beta>1$. 
The following theorem provides sufficient conditions under which the optimal $\Delta_i$ and $t_i$ can be uniquely found from the first-order conditions.
{Its proof can be found in
Appendix~\ref{app: proofs_wos}}.

\begin{theorem}\label{THM1}
Let ${\mathscr P}(\cdot)$ be strictly convex, and let $g(\Delta):={\mathscr P}'(\Delta)$ have a continuous inverse. Then \eqref{FOC1} and \eqref{FOC2} have a unique positive solution for $(t_i,\Delta_i)\in{\mathbb R}_+^2$,
{which yields the minimum ${\rm WOS}_n({\bs B},\omega,\mathscr{P})$.}
\end{theorem}

For instance, in the case of ${\mathscr P}(\Delta) = (\alpha/\beta) \Delta^\beta$ we obtain the first-order conditions
\[ (1-\omega)F_i(t_i) = \alpha\,\Delta_i^{\beta-1}, \quad \omega(1-F_i(t_i + \Delta_i)) = \alpha \,\Delta_i^{\beta-1}.\]
Solving $\Delta_i$ from the first of these two equations, we find
\[\Delta_i(t_i) = \left(\frac{1-\omega}{\alpha}F_i(t_i)\right)^{1/(\beta-1)},\]
which is an increasing function of $t_i$. Then $t_i$ is found from the equation ${\mathfrak f}^-_i(t_i)={\mathfrak f}^+_i(t_i)$, with ${\mathfrak f}^-_i(\cdot), {\mathfrak f}^+_i(\cdot)$ defined in Appendix~\ref{app: proofs_wos}.



\subsection{WOS under the normal distribution}\label{subsec:WOSnorm}
In the above, the travel times $B_i$ are defined in a general form. 
Here, we investigate the specific case that the travel-time distributions
are independent random variables with a normal distribution,
as, in that case, the distribution of $S_i$ is
normally distributed as well. 

Let $B_i \sim N(\mu_i,\sigma_i^2)$, 
with $B_1,B_2,\dots$ independent, and with $\mu_i>0$ and $\sigma_i\ll \mu_i$. The requirement $\sigma_i\ll \mu_i$ is intended to make sure that the probability mass on the negative numbers is effectively negligible; in practical terms $\sigma_i < \mu_i/3$ suffices, as this yields a probability of a negative travel time of just over $10^{-3}$. Considering the problem
${\rm WOS}_n({\bs B},\omega,{\mathscr P})$, the crucial observation is  that
$S_i\sim N(\tilde\mu_i,\tilde\sigma_i^2)$ with $\tilde{\mu}_i := \sum_{j = 1}^i \mu_j$ and 
$\tilde{\sigma}_i^2 := \sum_{j = 1}^i \sigma_j^2$, i.e.,
\begin{align*}
    F_i(t) = \Phi\left(\frac{t-\tilde{\mu}_i}{\tilde{\sigma}_i} \right),
\end{align*}
with $\Phi(\cdot)$ the cumulative distribution function of a standard normal random variable. 
From \eqref{FOC1} and \eqref{FOC2}, it follows that we are to solve the equations
\begin{align*}
    \omega\left(1-\Phi\left(\frac{t_i+\Delta_i-\tilde{\mu}_i}{\tilde{\sigma}_i} \right)\right)
    &= (1-\omega)\Phi\left(\frac{t_i-\tilde{\mu}_i}{\tilde{\sigma}_i} \right);\\
    \omega \left(1-\Phi\left(\frac{t_i+\Delta_i-\tilde{\mu}_i}{\tilde{\sigma}_i}\right)\right) &={\mathscr P}'(\Delta_i),
\end{align*}
if the optimum corresponds to a stationary point. 
To get a better intuition for these equations, consider again the specific example of ${\mathscr P}(\Delta) = \alpha \Delta$ for some $\alpha\leqslant \min\{\omega,1-\omega\}$.
We find that the window is
\[\left(\tilde\mu_i +\tilde\sigma_i\Phi^{-1}\left(\frac{\alpha}{1-\omega}\right),
\tilde\mu_i +\tilde\sigma_i\Phi^{-1}\left(1-\frac{\alpha}{\omega}\right)
\right);\]
evidently, if the left boundary of the interval happens to drop below 0, it is to be replaced by~$0$.
If $\omega\not=\frac{1}{2}$, the window is asymmetric around $\tilde\mu_i.$ Indeed, as is readily checked, the midpoint of the interval lies at 
\[\tilde\mu_i +\frac{\tilde\sigma_i}{2}\left(\Phi^{-1}\left(1-\frac{\alpha}{\omega}\right)-\Phi^{-1}\left(1-\frac{\alpha}{1-\omega}\right)\right).\]
This midpoint is larger than $\tilde\mu_i$ if $\omega<\frac{1}{2}$ and smaller than $\tilde\mu_i$ if $\omega>\frac{1}{2}$, in line with the weights assigned to the deliverer being late and early, respectively. Also observe that the width of the window is decreasing in $\alpha$: the more the width is penalized, the more narrow the optimal window is. In the special case that the service times are all $N(\mu,\sigma^2)$, then the window is of the form
\[\left(i\mu +\sigma\sqrt{i}\,\Phi^{-1}\left(\frac{\alpha}{1-\omega}\right),
i\mu +\sigma\sqrt{i}\,\Phi^{-1}\left(1-\frac{\alpha}{\omega}\right)
\right)\]
for $i\in\{1,\ldots,n\}$; again, if the left boundary is smaller than 0, it is replaced by $0$. We conclude that in this homogeneous scenario the windows are intervals around $i\mu $ (not necessarily symmetric) with a width that is proportional to $\sqrt{i}$.

\begin{figure}
    \centering
    \scalebox{0.6}{\input{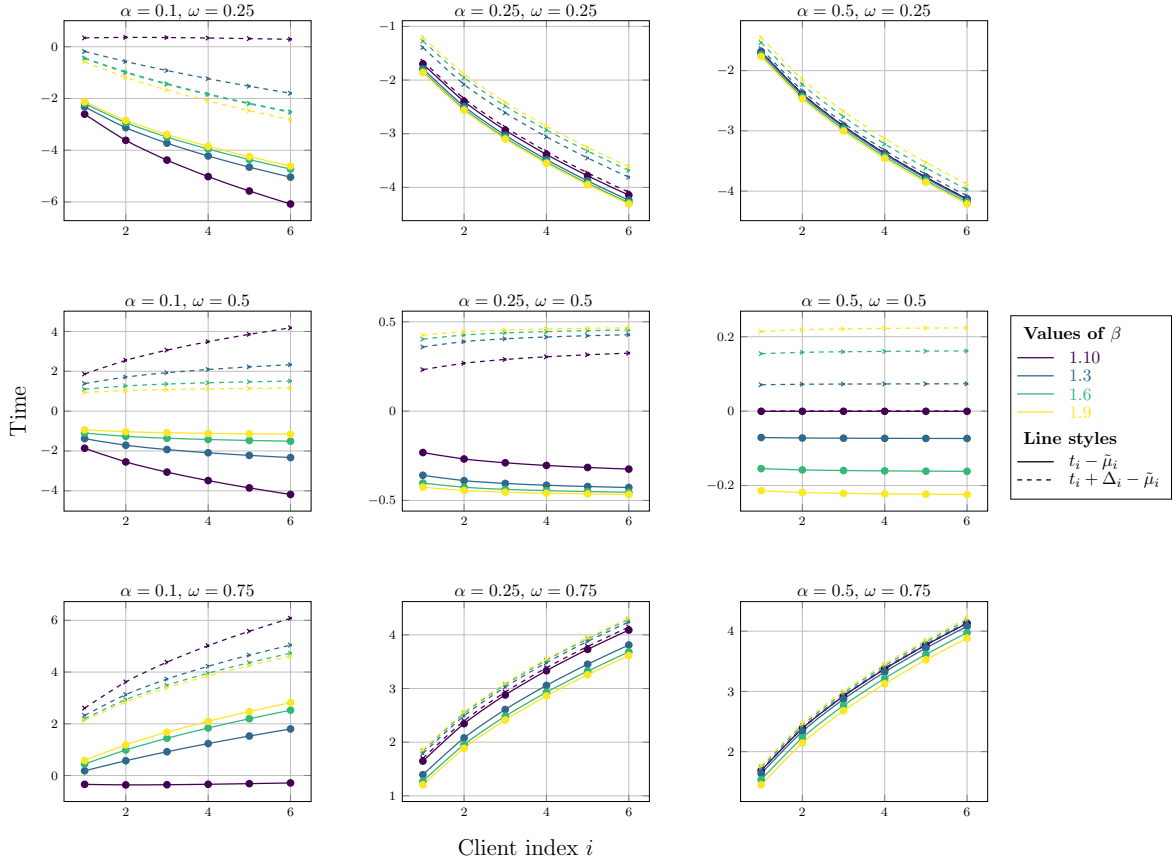}}
    \caption{Centered WOS time windows for normally distributed travel times with  parameters $\mu_i = 10$ and $\sigma_i = 2.5$. The three rows correspond to $\omega=0.25, 0.5, 0.75$, respectively, and the three columns to $\alpha = 0.1, 0.25, 0.5$, respectively. 
    \label{fig: normal_wos}}
\end{figure}

Whereas we find an explicit characterization of
the window in case $\mathscr{P}(\Delta) = \alpha \Delta$,
other penalty functions typically require numerical
optimization. Figure~\ref{fig: normal_wos} displays the `centered
windows' $[t_i - \tilde \mu_i, t_i + \Delta_i - \tilde \mu_i]$ 
for
various parameter values {in case ${\mathscr P}(\Delta) = (\alpha/\beta) \Delta^\beta$}. 
\begin{itemize}
    \item[$\circ$] 
It can be observed that the impact 
of $\omega$ aligns with the observations above:
$\omega = \frac{1}{2}$ yields a symmetric window,
but when $\omega < \frac{1}{2}$ (resp.\ $\omega > \frac{1}{2}$)
the penalty for early (resp.\ late) arrivals becomes dominant, 
and the optimal \(t_i\) shifts to an earlier 
(resp.\ later) time. 
We also observe symmetry
in the windows of $\omega$ and $1-\omega$, in that
the  centered window of $1-\omega$ is the mirrored
(in the time axis value of 0) centered window of $\omega$.
 \item[$\circ$] 
As \( \alpha \) scales the penalty linearly, higher values of $\alpha$
result in tighter time windows. 
The impact of $\beta$ is less trivial, as it is dependent on the value
of $\Delta$: 
the penalty for \( \Delta > 1 \) grows steeply with higher \( \beta \), 
but flattens for \( \Delta < 1 \).
Therefore, the order of the \( \Delta_i \) values with respect to \( \beta \) undergoes a transition, in that
for \( \alpha = 0.1 \) the smallest value of $\beta$ yields
the largest time-windows width, whereas
for $\alpha = 0.25, 0.5$ the smallest value of $\beta$ yields the smallest time-window width.
\end{itemize}


It is interesting to note that the symmetry of the windows around $\tilde{\mu}_i$ that is observed for $\omega = \frac{1}{2}$ holds more generally: it is true for any penalty function $\mathscr{P}(\cdot)$
satisfying the conditions of Theorem~\ref{THM1} and any service time
distribution that is symmetric in some $\tilde \mu_i$. 
As it turns out, the same conditions are also sufficient
for the mirrored windows of $\omega$ and $1-\omega$. The following lemma,
whose proof can be found in 
Appendix~\ref{app: proofs_wos},
makes these claims formal. 

\begin{lemma}\label{lem: symmetry}
Let ${\mathscr P}(\cdot)$ be strictly convex, 
$g(\Delta):={\mathscr P}'(\Delta)$ have a continuous inverse, 
and let the probability density function of $S_i$ be symmetric in $\tilde \mu_i \in \mathbb{R}_{+}$.
Then, if $[t(\omega'), t(\omega') + \Delta]$ is the
optimal window for $\omega = \omega'$, 
$[t(1-\omega'), t(1-\omega') + \Delta]$ is the
optimal window for $\omega = 1-\omega'$, where $t(1-\omega')$ is such that
\begin{align} \label{eq: omega_prime}
t(\omega') - \tilde \mu_i + \Delta = -(t(1-\omega') - \tilde \mu_i).
\end{align}
In the special case $\omega = \frac{1}{2}$, the window of client~$i$ is symmetric around
$\tilde \mu_i$.
\end{lemma}

\subsection{WOS under general convolutions}\label{subsec:WOSgen}
If a closed-form expression for the distribution function $F_i: \mathbb{R}_{+} \to [0,1]$ of $S_i$ is not available, it can be evaluated numerically by iteratively convolving discretized versions of the travel-time distributions. While this is in principle also possible when the travel times are dependent, in this section we focus on the case that they are independent. 
The idea is to simply discretize the continuous random variables $B_i$, $i = 1,\dots,n$, over a grid from $\mu_i-a_i$ to $\mu_i+a_i$,
for some $a_i \in \mathbb{R}_+$, using
a step size of $\varepsilon$. 
Then the arrival time distributions at the different client locations follow from a sequential numerical convolution
of these discrete distributions. Concretely, 
the arrival time distribution of client~$j$
follows from convolving the arrival time distribution
of client $j-1$ and the distribution $B_j$.
In all experiments below,
we have used $a_i = 4\sigma_i$ and a step size of $\varepsilon = 10^{-3}$.





\begin{figure}
    \centering
    \scalebox{0.7}{\begin{tikzpicture}
\begin{groupplot}[
    group style={group size=3 by 1, horizontal sep=2cm},
    width=7cm,
    height=8cm,
    xlabel={Client index $i$},
    ylabel={Time},
    xtick={1,2,3,4,5,6},
    grid=both,
    legend style={at={(1.9,1.13)}, anchor=south, legend columns=6, font=\small},
    title style={font=\small},
]

\nextgroupplot[title={$\omega=0.25$}]
\addplot+[thick, mark=o, blue]
    table [x=client, y=t_i, col sep=comma] {Figures_latex/comparison_plot_omega025_normal.csv};
\addlegendentry{Normal $t_i$}
\addplot+[thick, mark=x, blue, dashed]
    table [x=client, y=t_i_plus_delta, col sep=comma] {Figures_latex/comparison_plot_omega025_normal.csv};
\addlegendentry{Normal $t_i+\Delta_i$}
\addplot+[thick, mark=square*, orange]
    table [x=client, y=t_i, col sep=comma] {Figures_latex/comparison_plot_omega025_weibull.csv};
\addlegendentry{Weibull $t_i$}
\addplot+[thick, mark=star, orange, dashed]
    table [x=client, y=t_i_plus_delta, col sep=comma] {Figures_latex/comparison_plot_omega025_weibull.csv};
\addlegendentry{Weibull $t_i+\Delta_i$}
\addplot+[thick, mark=triangle*, green!70!black]
    table [x=client, y=t_i, col sep=comma] {Figures_latex/comparison_plot_omega025_lognormal.csv};
\addlegendentry{Lognormal $t_i$}
\addplot+[thick, mark=diamond*, green!70!black, dashed]
    table [x=client, y=t_i_plus_delta, col sep=comma] {Figures_latex/comparison_plot_omega025_lognormal.csv};
\addlegendentry{Lognormal $t_i+\Delta_i$}

\nextgroupplot[title={$\omega=0.5$}]
\addplot+[thick, mark=o, blue]
    table [x=client, y=t_i, col sep=comma] {Figures_latex/comparison_plot_omega05_normal.csv};
\addplot+[thick, mark=x, blue, dashed]
    table [x=client, y=t_i_plus_delta, col sep=comma] {Figures_latex/comparison_plot_omega05_normal.csv};
\addplot+[thick, mark=square*, orange]
    table [x=client, y=t_i, col sep=comma] {Figures_latex/comparison_plot_omega05_weibull.csv};
\addplot+[thick, mark=star, orange, dashed]
    table [x=client, y=t_i_plus_delta, col sep=comma] {Figures_latex/comparison_plot_omega05_weibull.csv};
\addplot+[thick, mark=triangle*, green!70!black]
    table [x=client, y=t_i, col sep=comma] {Figures_latex/comparison_plot_omega05_lognormal.csv};
\addplot+[thick, mark=diamond*, green!70!black, dashed]
    table [x=client, y=t_i_plus_delta, col sep=comma] {Figures_latex/comparison_plot_omega05_lognormal.csv};

\nextgroupplot[title={$\omega=0.75$}]
\addplot+[thick, mark=o, blue]
    table [x=client, y=t_i, col sep=comma] {Figures_latex/comparison_plot_omega075_normal.csv};
\addplot+[thick, mark=x, blue, dashed]
    table [x=client, y=t_i_plus_delta, col sep=comma] {Figures_latex/comparison_plot_omega075_normal.csv};
\addplot+[thick, mark=square*, orange]
    table [x=client, y=t_i, col sep=comma] {Figures_latex/comparison_plot_omega075_weibull.csv};
\addplot+[thick, mark=star, orange, dashed]
    table [x=client, y=t_i_plus_delta, col sep=comma] {Figures_latex/comparison_plot_omega075_weibull.csv};
\addplot+[thick, mark=triangle*, green!70!black]
    table [x=client, y=t_i, col sep=comma] {Figures_latex/comparison_plot_omega075_lognormal.csv};
\addplot+[thick, mark=diamond*, green!70!black, dashed]
    table [x=client, y=t_i_plus_delta, col sep=comma] {Figures_latex/comparison_plot_omega075_lognormal.csv};

\end{groupplot}
\end{tikzpicture}}
    \caption{
    {Centered WOS time windows using numerical convolutions; we use the same mean ($\mu_i = 10$) and standard deviation ($\sigma_i = 2.5$) for all distribution functions. }} 
    \label{fig:exp1_v1}
\end{figure}
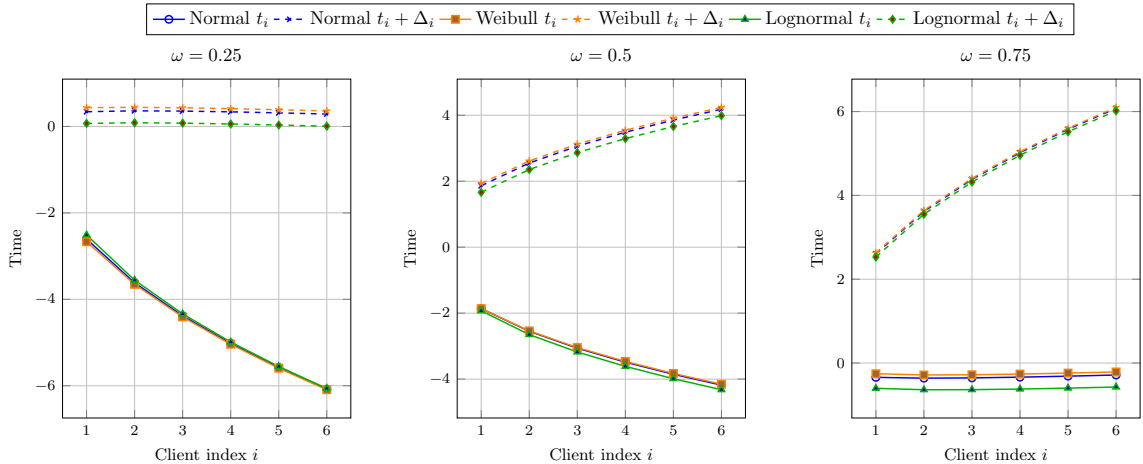

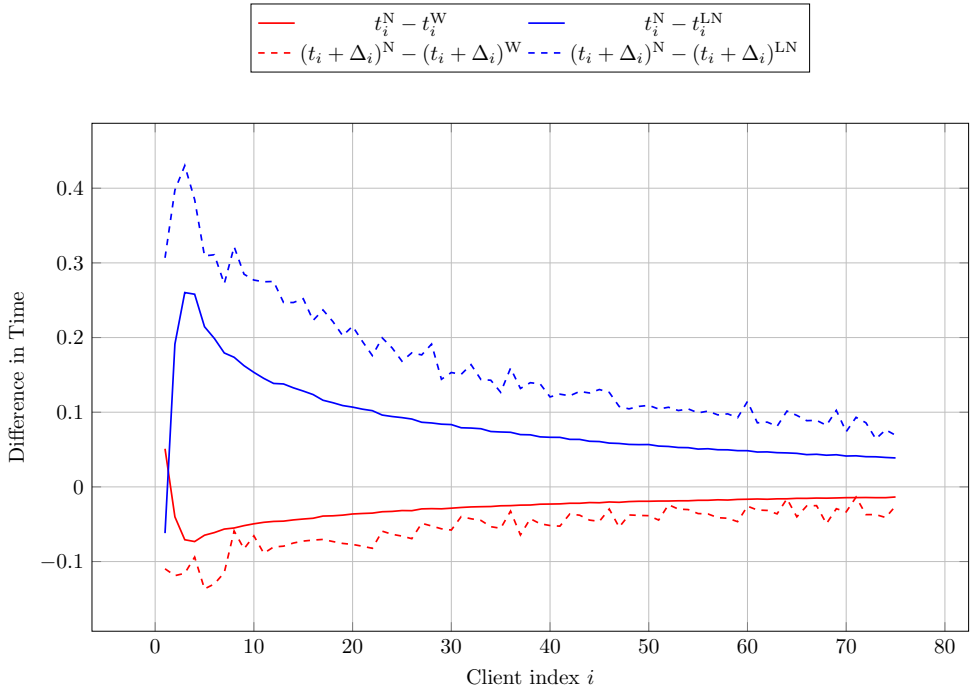
\begin{figure}
    \centering
    \scalebox{0.8}{\begin{tikzpicture}
\begin{axis}[
    width=16cm, height=10cm,
    xlabel={Client index $i$},
    ylabel={Difference in Time},
    grid=both,
    legend style={at={(0.5,1.10)}, anchor=south, font=\small, draw=black},
    legend columns=2,
    tick label style={font=\small},
    label style={font=\small},
    title style={font=\normalsize}
]
\addplot[
    red, thick
] table [
    x=client_index, y=diff_t_weibull, col sep=comma
] {Figures_latex/differences_comparison_new.csv};
\addlegendentry{$t_i^{\mathrm{N}} - t_i^{\mathrm{W}}$}

\addplot[
    blue, thick
] table [
    x=client_index, y=diff_t_lognormal, col sep=comma
] {Figures_latex/differences_comparison_new.csv};
\addlegendentry{$t_i^{\mathrm{N}} - t_i^{\mathrm{LN}}$}

\addplot[
    red, dashed, thick
] table [
    x=client_index, y=diff_window_weibull, col sep=comma
] {Figures_latex/differences_comparison_new.csv};
\addlegendentry{$(t_i+\Delta_i)^{\mathrm{N}} - (t_i+\Delta_i)^{\mathrm{W}}$}

\addplot[
    blue, dashed, thick
] table [
    x=client_index, y=diff_window_lognormal, col sep=comma
] {Figures_latex/differences_comparison_new.csv};
\addlegendentry{$(t_i+\Delta_i)^{\mathrm{N}} - (t_i+\Delta_i)^{\mathrm{LN}}$}

\end{axis}
\end{tikzpicture}}
    \caption{Difference of time windows for Weibull and Lognormal distribution, relative to time windows for normal distribution, for $\omega = 0.25$. (i)~Red: the $t_i$ based on the normal distribution minus those based on the Weibull distribution, (ii)~red dotted: same, but then for the $t_i+\Delta_i$, (iii)~blue: the $t_i$ based on the normal distribution minus those based on the Lognormal distribution, (iv)~blue dotted: same, but then for the $t_i+\Delta_i$.}
    \label{fig:normal_approx_difference}
\end{figure}


Figure~\ref{fig:exp1_v1} illustrates the time windows resulting from the convolution method across three different distributions: normal, Weibull, and Lognormal. As before, these time windows have been `centered', meaning that the mean $\tilde \mu_i$ has been subtracted for each client $i$. {The mean and standard deviation of the service times are 10 and 2.5, respectively, for all clients and distributions. With these parameters, the normal and Weibull distributions exhibit similar, symmetric shapes, resulting in nearly identical time windows.}
In contrast, the density of the Lognormal distribution is significantly skewed, leading to {somewhat different} 
time windows. Specifically, for $\omega = \frac{1}{2}$, the Lognormal time windows are not symmetric around zero.
It is worth noting that, for all the settings shown in Figure~\ref{fig:exp1_v1}, the endpoints of the time windows under the Lognormal distribution occur earlier than those under the Weibull and normal distributions. This effect arises because the Lognormal distribution has a greater probability mass to the left of the mean, compared to the other two distributions.

Figure~\ref{fig:exp1_v1} considers the case of $n = 6$ customers. 
{Although computations remain feasible for larger values of $n$, the numerical evaluation of convolutions becomes increasingly demanding, as the number of bins required grows with the number of clients.} 
One potential speed-up, at the cost of a small loss
in accuracy, is to bound the total number of bins that may arise after a convolution step by changing the bin-size and allowing for larger bins, i.e., by increasing $\varepsilon$, thereby re-distributing the weights of the bins (e.g., 
using interpolation). 

A more convenient solution to cope with the computational burden, however, is to approximate, for sufficiently high client indices $i$, the $S_i$ by their normally distributed counterparts. 
This approximation is motivated by the central limit theorem, noting that $S_i$ is the sum of random variables. 
In particular, for any client $i\geqslant i_0$, 
we set $S_i \sim N(\tilde {\mu}_i, {\tilde \sigma}_i^2)$; in the literature often $i_0=15$ is being used. The idea is that the use of this normal approximation, {which
we will refer to as the
$N(i_0)$-approximation},
significantly reduces computational costs
while maintaining reasonable accuracy.

The accuracy of the $N(15)$-approximation is illustrated in Figure~\ref{fig:normal_approx_difference} for cases in which the $B_i$ follow a Weibull or Lognormal distribution.
The figure shows a typical example of the differences in the resulting time windows.
A comparison in the value of the objective function (based on $R=10^9$ simulation runs)
for the same distributions is given in Table~\ref{tab:performance_runtime_2}.
Specifically, for both distributions, time windows are computed using {(a)~the convolution method, 
(b)~the $N(1)$-approximation,
and (c)~the $N(15)$-approximation.}
For these sets of time windows, we evaluate the value of the objective function $G_n(\boldsymbol{t}, \boldsymbol{\Delta})$ using Monte Carlo simulation, that is, by $R$ times sampling the $n$ 
random variables $B_i$, evaluating for those sampled $B_i$ the value of the objective function, and taking their average (that is, add the $R$ realized values of the cost function and divide by $R$). {We also calculated the percentage error of the $N(1)$-
and $N(15)$-approximation relative to the (exact) convolution-based value, respectively indicated by $\Delta_{N(1)}^i$ and $\Delta_{N(15)}^i$;} here $i\in\{{\rm W},{\rm LN}\}$, where W stands for Weibull and LN for Lognormal. 
As can be observed from Table~\ref{tab:performance_runtime_2},
for both $B_i$ having a Lognormal and a Weibull distribution,
there are virtually negligible differences between the output obtained via the three methods. 
{Considering that the convolution method is significantly more time-consuming than the normal approximations, we opt to use the normal approximations.}

In Table~\ref{tab:performance_runtime_2}, the travel times are assumed i.i.d. To stress-test the normal approximations, Table~\ref{tab:non_iid_variance_sin} considers a setting where the variances of the $B_i$'s vary with the client index, where its standard deviation is given by $2.5 + 1.5 \sin\left({2\pi i}/{10}\right)$. Even in this heterogeneous setting, the $N(15)$-approximation remains highly accurate. The run times are omitted because they are virtually identical to those reported in Table~\ref{tab:performance_runtime_2}.

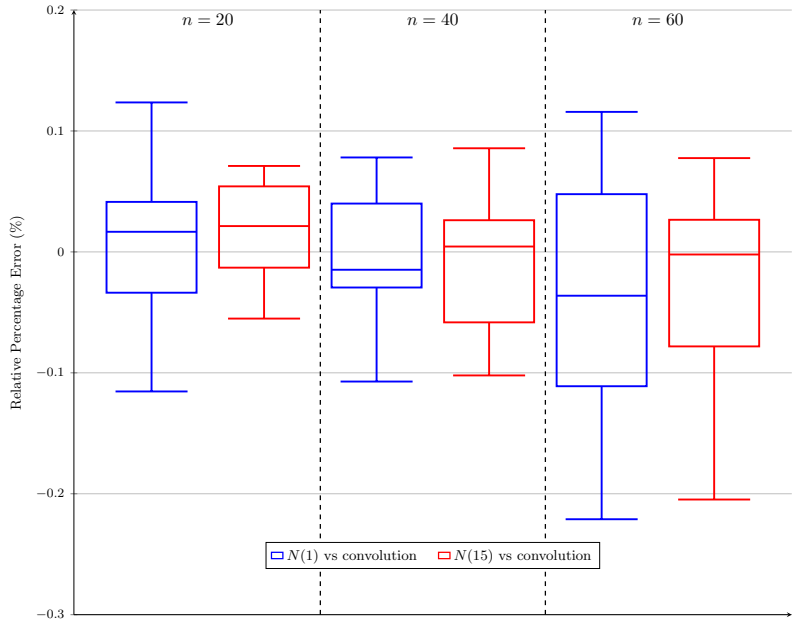
\begin{figure}
    \centering
    \scalebox{0.6}{\begin{tikzpicture}
\begin{axis}[
  boxplot/draw direction=y,
  axis x line=bottom, axis y line=left,
  ylabel={Relative Percentage Error (\%)},
  ymin=-0.30, ymax=0.20,          
  width=1.05\textwidth,
  enlarge x limits=0.05,
  clip mode=individual,            
  ymajorgrids,
  xtick=\empty,
  ytick={-0.3,-0.2,-0.1,0,0.1,0.2},
  legend columns=2,
  legend style={at={(0.5,0.12)}, anchor=north, /tikz/every even column/.append style={column sep=1em}},
  legend image code/.code={\draw[#1,thick] (0cm,-0.07cm) rectangle (0.25cm,0.07cm);}
]

\addplot+[
  boxplot prepared={
    median=0.016624,
    upper quartile=0.041348,
    lower quartile=-0.033751,
    upper whisker=0.123597,
    lower whisker=-0.115403
  },
  draw=blue
] coordinates {};
\addlegendentry{$N(1)$ vs convolution}

\addplot+[
  boxplot prepared={
    median=0.021347,
    upper quartile=0.054184,
    lower quartile=-0.013034,
    upper whisker=0.071066,
    lower whisker=-0.055148
  },
  draw=red
] coordinates {};
\addlegendentry{$N(15)$ vs convolution}

\addplot+[
  boxplot prepared={
    median=-0.014768,
    upper quartile=0.039965,
    lower quartile=-0.029449,
    upper whisker=0.078093,
    lower whisker=-0.107209
  },
  draw=blue
] coordinates {};

\addplot+[
  boxplot prepared={
    median=0.004439,
    upper quartile=0.026193,
    lower quartile=-0.058269,
    upper whisker=0.085692,
    lower whisker=-0.102153
  },
  draw=red
] coordinates {};

\addplot+[
  boxplot prepared={
    median=-0.036237,
    upper quartile=0.047738,
    lower quartile=-0.111085,
    upper whisker=0.115782,
    lower whisker=-0.221075
  },
  draw=blue
] coordinates {};

\addplot[
  boxplot prepared={
    median=-0.002149,
    upper quartile=0.026553,
    lower quartile=-0.078152,
    upper whisker=0.077564,
    lower whisker=-0.204774
  },
  draw=red
] coordinates {};


\draw[dashed, thick] (axis cs:2.5,\pgfkeysvalueof{/pgfplots/ymin})
                     -- (axis cs:2.5,\pgfkeysvalueof{/pgfplots/ymax});
\draw[dashed, thick] (axis cs:4.5,\pgfkeysvalueof{/pgfplots/ymin})
                     -- (axis cs:4.5,\pgfkeysvalueof{/pgfplots/ymax});

\node at (axis cs:1.5, \pgfkeysvalueof{/pgfplots/ymax}) [yshift=-6pt] {\textbf{$n=20$}};
\node at (axis cs:3.5, \pgfkeysvalueof{/pgfplots/ymax}) [yshift=-6pt] {\textbf{$n=40$}};
\node at (axis cs:5.5, \pgfkeysvalueof{/pgfplots/ymax}) [yshift=-6pt] {\textbf{$n=60$}};

\end{axis}
\end{tikzpicture}}
    
     \caption{ Boxplot assessing accuracy of the normal approximation over a broad range of heterogeneous scenarios. The hinges of the box correspond to the 25- and 75-percentiles, the whiskers to the 5- and 95-percentiles. 
   The predefined parameter sets used in experiments are: $\omega\in \{0.25, 0.5, 0.75\}$, $\alpha\in\{0.01, 0.05, 0.1, 0.15\}$,   $\beta\in\{1.1, 1.2, 1.3, 1.4, 1.5\} $, $\mu\in\{10, 12, 15\}$.}
    \label{fig:boxplot_errors}
\end{figure}

In the above experiments, we picked specific values for the parameters of the underlying model. The last experiment of this section aims to show that the above conclusions carry over to a broad range of such parameters. 
Figure~\ref{fig:boxplot_errors} describes the effect of stronger client heterogeneity,
for three values of the number of clients $n$ and estimates based on simulation runs $R = 10^7$. In each run, the values of $\omega, \alpha$, and $\beta$ are sampled uniformly at random from predefined sets (given in the caption of the figure).
Likewise, each travel-time distribution is
randomly selected to be normal, Lognormal or Weibull, with mean $\mu$
sampled, again uniformly at random, from predefined sets as well (also given in the caption of the figure); the standard deviation is kept at 2.5 for each travel time. 

For each experiment, we compute the objective function, with (a)~the windows based on the normal approximations and  (b)~the windows based on the convolution method, and then we compute their percentage difference,
so as to assess the accuracy of the $N(1)$- and $N(15)$-approximations. These percentage differences are shown in the boxplot of Figure~\ref{fig:boxplot_errors}. 
The main conclusion from the boxplot is that both the $N(1)$- and $N(15)$-approximations are highly accurate across the broad range of heterogeneous scenarios that were studied.

{\small
\begin{table}
    \centering
    { 
    \begin{tabular}{lrrrrr}
        \toprule
        \multicolumn{1}{c}{\textbf{}} & \multicolumn{5}{c}{\textit{Objective function}} \\ 
        \cmidrule(lr){2-6} 
        number of clients $n$& 20 & 40 & 60 & 80 & 100  \\ 
         \midrule
        {Weibull -- convolution} & 54.4714 & 155.8506 & 288.7246 & 447.1245 & 627.7793  \\
        & \scriptsize{(11.46)} & \scriptsize{(39.13)} & \scriptsize{(89.00)} & \scriptsize{(138.91)} & \scriptsize{(215.22)}  \\
        \midrule
        {Weibull -- $N(1)$-approximation} & 54.3346 & 155.7097 & 288.5674 & 447.0486 & 627.6036  \\
        &\scriptsize{(0.06)} & \scriptsize{(0.07)} & \scriptsize{(0.10)} & \scriptsize{(0.14)} & \scriptsize{(0.16)}  \\
        \midrule
        {Weibull -- $N(15)$-approximation} & 54.4300 & 155.7639 & 288.3891 & 447.0403 & 627.6280  \\
         & \scriptsize{(6.12)} & \scriptsize{(6.24)} & \scriptsize{(6.67)} & \scriptsize{(6.68)} &\scriptsize{ (6.70)}  \\
         \midrule
        {Lognormal -- convolution} & 54.2104 & 155.4583 & 288.2273 & 446.6121 & 627.2115   \\
         & \scriptsize{(10.19)} & \scriptsize{(37.03)} & \scriptsize{(80.23)} & \scriptsize{(136.57)} & \scriptsize{(207.60)}   \\
        \midrule
        {Lognormal -- $N(1)$-approximation} & 54.4198& 155.7562& 288.4493& 446.4572& 627.8023  \\
         &\scriptsize{(0.06)} & \scriptsize{(0.06)} &\scriptsize{(0.08)}& \scriptsize{(0.11)}&\scriptsize{(0.14)}   \\
        \midrule
        {Lognormal -- $N(15)$-approximation} & 54.2375 & 155.5400 & 288.4854 & 447.0526 & 626.8255   \\
         & \scriptsize{(6.08)} & \scriptsize{(6.09)} & \scriptsize{(6.14)} & \scriptsize{(6.44)} & \scriptsize{(6.44)}   \\
        \midrule 
       \textbf{$\Delta_{ N(1)}^{\rm W}$}  &-0.25	&-0.09	&-0.05	&-0.02	&-0.003	   \\
       \textbf{$\Delta_{ N(15)}^{\rm W}$} &-0.08 &-0.06 &-0.12 &-0.02 &-0.02   \\
        \textbf{$\Delta_{N(1)}^{\rm LN}$} & 0.38& 0.19&0.08 & -0.03& 0.09  \\
        \textbf{$\Delta_{N(15)}^{\rm LN}$} &0.04 &0.05 & 0.09 & 0.10 &-0.06    \\
        \bottomrule
    \end{tabular}}
    \caption{Accuracy  comparisons, for $B_i$\,s i.i.d., with $\mu_i=10$ and  $\sigma_i=2.5$, $\alpha = 0.1$, $\beta = 1.5$; runtimes (in seconds) have been added between brackets, below the respective objective function values.}
    \label{tab:performance_runtime_2}
\end{table}}


{\small
\begin{table}[htbp]
\centering
\caption{Accuracy comparisons, for $B_i$\,s non-i.i.d., with $\mu_i=10$ and  $\sigma_i=2.5 + 1.5 \sin\left(\frac{2\pi i}{10}\right)$, $\alpha = 0.1$, $\beta = 1.5$.}
{%
\begin{tabular}{lccccc}
\toprule
 & \multicolumn{5}{c}{\it Objective function} \\ 
\cmidrule(lr){2-6}
number of clients $n$ & 20 & 40 & 60 & 80 & 100 \\
\midrule
Weibull -- convolution     & 40.6738 & 117.8401 & 219.3660 & 340.7553 & 479.2350 \\
Weibull -- $N(1)$-approximation                & 40.5520 & 117.6964 & 219.1714 & 340.5374 & 479.4331 \\
Weibull -- $N(15)$-approximation  & 40.6108 & 117.5708 & 219.0133 & 340.5200 & 478.4872 \\
\midrule 
Lognormal -- convolution   & 40.5212 & 117.5882 & 219.0586 & 340.4068 & 478.8196 \\
Lognormal -- $N(1)$-approximation                 & 40.5922 & 117.6391 & 219.1380 & 340.7982 & 478.9577 \\ 
Lognormal -- $N(15)$-approximation  & 40.5094 & 117.5752 & 219.1287 & 340.6473 & 478.5470  \\
\bottomrule
\end{tabular}%
}
\label{tab:non_iid_variance_sin}
\end{table}}


\section{Optimal static windows: the UWOS case}\label{sec:uwos}
In this section, we numerically evaluate window assignments based on the UWOS optimization problem. Some delivery companies assign equally wide time windows to all clients in an effort to ensure fairness. While this approach may promote fairness in terms of window width, it is important to note that other factors -- specifically, the deliverer's arrival before or after the announced window -- can vary significantly between clients. Here, we quantify the degree of fairness (or unfairness) experienced by different clients.

It is readily checked that in this case with uniform window lengths, the first-order conditions are given by the following $n+1$ conditions:
\begin{align*}
    \frac{\partial}{\partial t_i} G_n({\bs t},\Delta {\bs 1}) &=  (1-\omega)F_i(t_i)- \omega(1-F_i(t_i + \Delta))=0 ,\\
    \frac{\partial}{\partial \Delta} G_n({\bs t},{\Delta}{\bs 1}) &= - \omega\sum_{i=1}^n(1-F_i(t_i + \Delta)) + n{\mathscr P}'(\Delta)=0.
\end{align*}
As in the WOS case, we start our analysis with a linear window cost, i.e., ${\mathscr P}(\Delta)=\alpha\Delta$. These equations allow for a sequential solution. To this end, we define
\[H_{\bs t}(\Delta):=\sum_{i=1}^n F_i(t_i + \Delta).\]
Upon combining the first-order equations, we find that the optimal width for a given vector ${\bs t}$ is 
\begin{align} \label{eq: delta_uwos_linear}
\Delta({\bs t}) =  H_{\bs t}^{-1}\left(\frac{n(\omega-\alpha)}{\omega}\right).
\end{align}
The optimal vector ${\bs t}$ follows from the $n$ coupled equations
\[(1-\omega)F_i(t_i)= \omega(1-F_i(t_i + \Delta({\bs t}))).\]
Similar to the WOS objective, $\alpha \leqslant \min\{\omega, 1-\omega\}$ is a necessary requirement for the existence of a solution. Indeed, since $H_{\bs t}: \mathbb{R}_{\geqslant 0} \to [0, 1]$, 
\eqref{eq: delta_uwos_linear} yields 
$\alpha \leqslant \omega$. Moreover, combining the $n+1$
conditions also gives $(1-\omega) H_{\bs t}(0) 
    = n\alpha$,
from which $\alpha \leqslant 1-\omega$ follows.

Note that if $\mathscr{P}'(\Delta)$ has a continuous inverse and $\mathscr{P}(\cdot)$ is strictly
convex, we can write
\begin{align*}
    \Delta({\bs t}) = g^{-1}\left(\frac{1-\omega}{n}
    \sum_{i = 1}^n F_i(t_i)\right).
\end{align*}
Then, the $n$ equations that need to
be solved are of the form
\begin{align*}
    (1-\omega)F_i(t_i) = \omega(1-F_i(t_i + \Delta(\bs t))).
\end{align*}
In order to (numerically) solve this system, a particularly practical property is that
in case all $t_j$ are fixed except $t_i$,
the same arguments as those used
in the proof of Theorem 1 establish the
existence of a unique solution $t_i$.

\begin{figure}
    \centering
    \scalebox{0.8}{\input{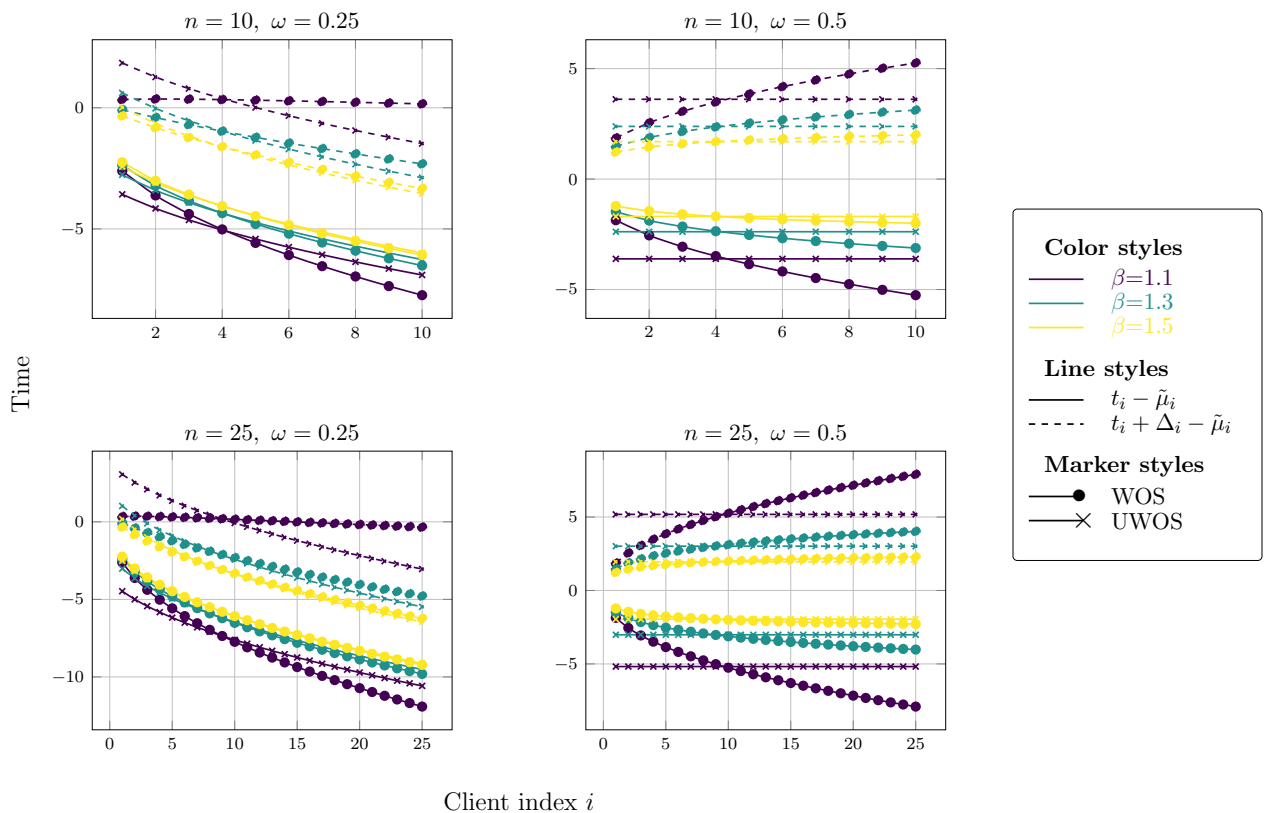}}
    \caption{Time windows UWOS vs.\ WOS. Per panel there are 12 lines: (a)~the bullets correspond to WOS, the crosses to UWOS; (b)~the solid lines correspond to the start of the windows, the dashed lines to the end of the windows; (c)~the three colors correspond to the three values of $\beta$.}
    \label{fig:time_windows_uwos_vs_wos}
\end{figure}

We proceed by discussing a number of insightful numerical experiments. In the first experiment, we present time windows as produced by UWOS with their counterparts produced by WOS. These are presented in Figure~\ref{fig:time_windows_uwos_vs_wos}, in which the travel times are i.i.d.\ with mean $\mu=10$ and standard deviation $\sigma=2.5$. The graphs confirm the conceivable properties that, across all scenarios considered, (i) for WOS we have that  $\Delta_i$
is a non-decreasing function, 
and (ii)~$\Delta_1 \leqslant \Delta \leqslant \Delta_n$, with $\Delta$ the window width for UWOS.

Since the parameter space optimized in WOS is larger than that in UWOS, it is natural to ask how much lower the objective value is under WOS 
compared to UWOS 
. Figure \ref{fig:uwos-wos-rel-perf-both} reveals that this difference is typically in the order of 10\%; the travel times are taken from the same distributions as in the experiments underlying Figure \ref{fig:time_windows_uwos_vs_wos}. 

\begin{figure}[t]
    \centering
    \scalebox{0.8}{\begin{tikzpicture}
\begin{groupplot}[
  group style={group size=3 by 1, horizontal sep=1.8cm},
  width=6.8cm, height=6.2cm,
  grid=both,
  xlabel={}, ylabel={},
  title style={yshift=-0.6em},
  table/col sep=comma
]

\def\fileTotal{Figures_latex/ti1/firstcol_total_n25_omega0_25.csv}
\def\fileStoch{Figures_latex/ti1/firstcol_stochastic_n25_omega0_25.csv}
\def\filePen{Figures_latex/ti1/firstcol_penalty_n25_omega0_25.csv}

\nextgroupplot[ylabel={Total objective}]
\plotPair{\fileTotal}{WOS_b11}{UWOS_b11}{bcolA1}{}
\plotPair{\fileTotal}{WOS_b13}{UWOS_b13}{bcolB1}{}
\plotPair{\fileTotal}{WOS_b15}{UWOS_b15}{bcolC1}{}

\nextgroupplot[ylabel={Stochastic term}, xlabel={Client index $i$}]
\plotPair{\fileStoch}{WOS_b11}{UWOS_b11}{bcolA1}{}
\plotPair{\fileStoch}{WOS_b13}{UWOS_b13}{bcolB1}{}
\plotPair{\fileStoch}{WOS_b15}{UWOS_b15}{bcolC1}{}

\nextgroupplot[ylabel={Penalty term}]
\plotPair{\filePen}{WOS_b11}{UWOS_b11}{bcolA1}{}
\plotPair{\filePen}{WOS_b13}{UWOS_b13}{bcolB1}{}
\plotPair{\filePen}{WOS_b15}{UWOS_b15}{bcolC1}{}

\end{groupplot}

\node[
  draw=black, fill=white, rounded corners=2pt,
  inner sep=3pt, anchor=south, scale=0.9
] at ($(group c2r1.north) + (0,0.7)$) {%
\begin{tikzpicture}[baseline]
  \draw[thick, color=bcolA1] (0,0) -- +(0.9,0); \node[anchor=west] at (1.1,0) {$\beta=1.1$};
  \draw[thick, color=bcolB1] (3.2,0) -- +(0.9,0); \node[anchor=west] at (4.4,0) {$\beta=1.3$};
  \draw[thick, color=bcolC1] (6.4,0) -- +(0.9,0); \node[anchor=west] at (7.6,0) {$\beta=1.5$};

  \draw[very thick, color=black] (9.8,0) -- +(0.9,0);
  \draw[black] plot[mark=*, mark options={solid}] coordinates { (10.25,0) };
  \node[anchor=west] at (11.0,0) {WOS};

  \draw[thick, dashed, color=black] (13.2,0) -- +(0.9,0);
  \draw[black] plot[mark=x, mark options={solid}] coordinates { (13.65,0) };
  \node[anchor=west] at (14.6,0) {UWOS};
  
\end{tikzpicture}
};

\end{tikzpicture}}
    \caption{Per-client contribution to the objective function for UWOS and WOS for $\omega=0.25$.}
    \label{fig:uwos-wos-perf_per_client}
\end{figure}
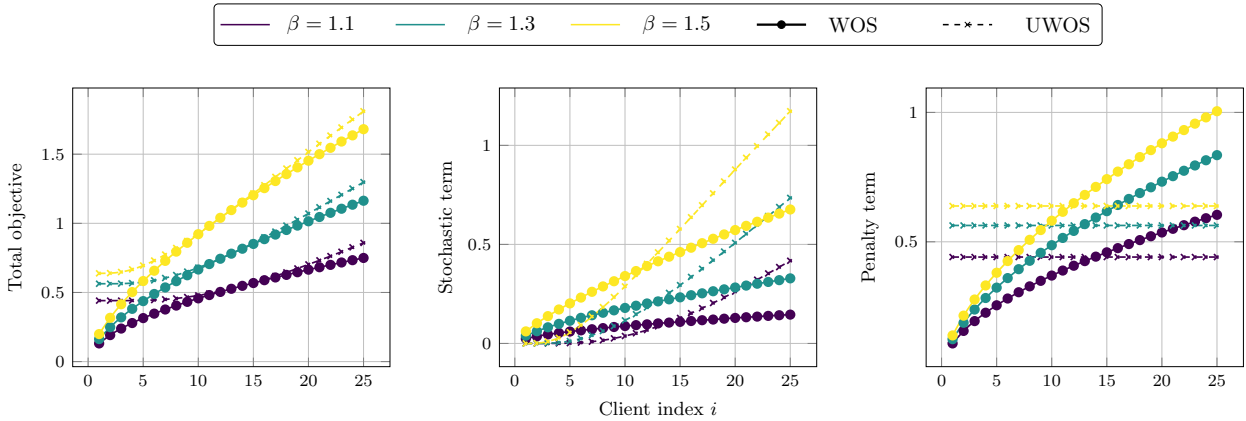

In UWOS, all clients are assigned time windows of equal length. An important question is to what extent this leads to fairness among clients, especially considering that the other contributions to the objective function, resulting from the deliverer arriving too early or too late, can still vary from one client to another. 
The right-most plot of Figure~\ref{fig:uwos-wos-perf_per_client} shows that UWOS is (evidently) more fair as it comes to the contribution to the window width. However, the left-most plot shows that UWOS has a lower per-client spread of the disutilities than WOS, but the clients early in the schedule are still considerably better of than the later clients.

\begin{figure}
    \centering
    \scalebox{0.8}{\begin{tikzpicture}
\begin{groupplot}[
  group style={group size=2 by 1, horizontal sep=2.2cm},
  width=7.5cm, height=5.5cm,
  grid=both,
  xlabel={Number of clients $n$},
  ylabel={Performance}, 
  title style={yshift=-0.6em},
  table/col sep=comma
]

\def\fileA{Figures_latex/ti1/perf_omega0_25.csv}
\def\fileB{Figures_latex/ti1/perf_omega0_5.csv}

\nextgroupplot[title={$\,\omega=0.25$}]
\plotAbsPair{\fileA}{WOS_b11}{UWOS_b11}{bcolA1}{1.8pt}
\plotAbsPair{\fileA}{WOS_b13}{UWOS_b13}{bcolB1}{1.8pt}
\plotAbsPair{\fileA}{WOS_b15}{UWOS_b15}{bcolC1}{1.8pt}

\nextgroupplot[title={$\,\omega=0.5$}]
\plotAbsPair{\fileB}{WOS_b11}{UWOS_b11}{bcolA1}{1.8pt}
\plotAbsPair{\fileB}{WOS_b13}{UWOS_b13}{bcolB1}{1.8pt}
\plotAbsPair{\fileB}{WOS_b15}{UWOS_b15}{bcolC1}{1.8pt}

\end{groupplot}

\node[
  draw=black, fill=white, rounded corners=2pt,
  inner sep=2pt, anchor=west, scale=0.85
] at ($(group c2r1.east) + (0.9,0)$) {%
{\scriptsize
\begin{tikzpicture}[baseline]
  \node[anchor=west] at (0,0) {\textbf{$\beta$ colors}};
  \draw[thick, color=bcolA1] (0,-0.5) -- +(0.8,0); \node[anchor=west, color=bcolA1] at (0.95,-0.5) {$1.1$};
  \draw[thick, color=bcolB1] (0,-0.9) -- +(0.8,0); \node[anchor=west, color=bcolB1] at (0.95,-0.9) {$1.3$};
  \draw[thick, color=bcolC1] (0,-1.3) -- +(0.8,0); \node[anchor=west, color=bcolC1] at (0.95,-1.3) {$1.5$};
  \node[anchor=west] at (0,-2.0) {\textbf{Schemes}};
  \draw[thick, color=black] (0,-2.5) -- +(0.8,0); \node[anchor=west] at (0.95,-2.5) {WOS};
  \draw[thick, dashed, color=black] (0,-2.9) -- +(0.8,0); \node[anchor=west] at (0.95,-2.9) {UWOS};
\end{tikzpicture}
}
};

\end{tikzpicture}}
    \caption{Performance of UWOS and WOS for different
    values of $n$. The bullets correspond to WOS, the crosses to UWOS.}
    \label{fig:uwos-wos-rel-perf-both}
\end{figure}
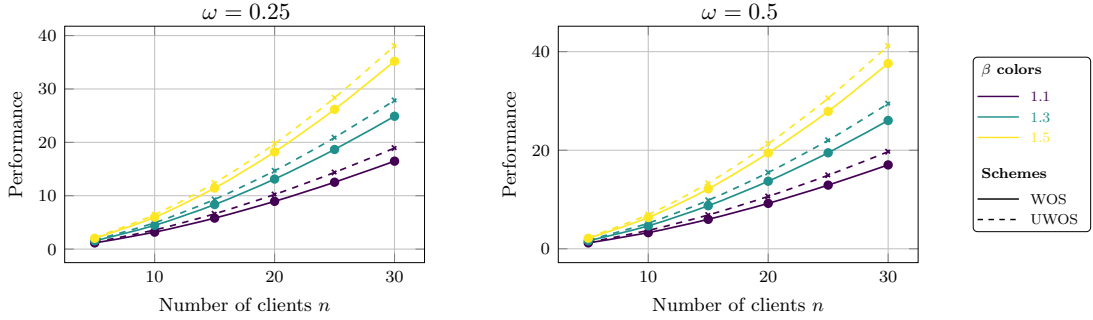






\section{Optimal dynamic windows: the DWOS case}\label{sec:dyn}

In this section we study the dynamic procedure that has been verbally introduced in Section~\ref{subsec:dyn}. 
In the computations presented below we work in the setting that all travel times are independent normally distributed random variables, which we have shown to be an accurate proxy in Section~\ref{subsec:WOSgen}, but in principle one could also work with generally distributed travel times using the machinery of Section~\ref{subsec:WOSgen}.

Recall that the idea of a dynamic procedure is that by the information that becomes available during the delivery process, the driver may
provide clients with updated time windows that reflect this information.
More specifically, if the driver is currently on its way to client $N$,
the travel times $B_1,\dots,B_{N-1}$ are known, and hence $S_{N-1}$ as well, so that there is less variability in the arrival time for clients $ N,\ldots,n$.
Hence, if the driver would compute
the optimal time windows at this time
instant, these would typically be more narrow than those computed at $t = 0$.
Moreover, if the realization of $S_{N-1}$ is right (resp.\ left) from its
mean, a newly computed time window would start later (resp.\ earlier).

Whereas the driver could, in principle,
send frequent updates of the time windows to the clients, having too many communication moments could be disruptive for the clients. Moreover,
communicating an updated time window at a late stage may be frustrating
for the clients; e.g., in case the updated time windows start earlier, 
it does not allow them to travel home in time. As discussed in Section~\ref{subsec:dyn}, if time windows can be updated during the delivery process while accounting for communication costs, dynamic programming could, in principle, be used to compute the optimal time window assignment at any given moment. This would rely on knowing the index of the current customer being `served' and the corresponding elapsed travel time. However, the dimensionality of the resulting optimization problem is prohibitively high, making this approach impractical in real-world operational settings.

To use the updated information about the executed part of the trip while
providing little and timely updates,  we propose DWOS: to use the static procedure dynamically, but only communicate an update the first time
the difference between the start of 
the updated window and the current time is smaller than some threshold $T$. More specifically, the efficiency
of the static procedure (especially with the use of the $N(1)$-approximation) allows us to recompute 
the time windows every $\tau$ time units. Define by $N_j$ the index of the customer in service by time $\tau j$, i.e., at the $j$-th recomputation:
\[N_j:=\max\left\{i:\sum_{k=1}^i B_k\leqslant \tau j\right\}+1.\]
At this recomputation epoch, we are given the {\it elapsed} travel time towards customer $N_j$; say it has the value \[b_j:=\tau j- \sum_{k=1}^{N_j-1}B_k.\] Then it follows directly that her residual travel time is characterized via the complementary distribution function, for $t\geqslant 0$,
\begin{equation}\label{eq:conditional}{\mathbb P}(B_{N_j}^{({\rm res})}>t) ={\mathbb P}(B_{N_j}>t+b_j\,|\,B_{N_j}>b_j) = \frac{{\mathbb P}(B_{N_j}>t+b_j)}{{\mathbb P(B_{N_j}>b_j)}};\end{equation}
the corresponding mean and variance can be evaluated by applying the expressions in Appendix \ref{app:cond}.
These allow for the $N(1)$-approximation of the arrival time distribution for clients $i \geq N_j$.
Let $(t_i, t_i + \Delta_i)$ be the 
time interval that minimizes WOS for
client~$i$ at time $\tau j$. Then, DWOS is such that 
client~$i$ receives an updated travel time if $t_i - \tau j \leq T$ and
if they have not received an update before.

The algorithm below formalizes the DWOS procedure. We use the notation 
\[{\bs B}[N]_j:=(B_{N_j}^{({\rm res})}, B_{N_j+1},\ldots,B_n),\]
which is a vector of dimension $M_j:=n -N_j+1$.  The dummy variable $\texttt{Flag}_i$, for $i=1,\ldots,n$, indicates whether client $i$ has been sent a time window update.

\begin{algorithm}[h]
\KwResult{Time window $(t_i, t_i + \Delta_i)$ for each client $i = 1,\dots,n$}

\textbf{a. Initialization:} Set $\texttt{Flag}_i = 0$ for all $i = 1,\dots,n$; set $N_0 = 1$, $j = 0$, and $b_0 = 0$\;

\textbf{b.} Let $(\boldsymbol{t}, \boldsymbol{\Delta}) \in \mathscr{U}_n$ be a minimizer of $\text{WOS}_n(\boldsymbol{B}, \omega, \mathscr{P})$\;

\While{\tt{\tt Flag}$_n = 0$}{
    \textbf{1.} Update $N_j$ and $b_j$\;

    \textbf{2.} Let $(\boldsymbol{t}, \boldsymbol{\Delta}) \in \mathscr{U}_{n - N_j + 1}$ be a minimizer of $\text{WOS}_{n - N_j + 1}(\boldsymbol{B}[N]_j, \omega, \mathscr{P})$\;

    \textbf{3.} \For{$i = N_j,\dots,n$}{
        \If{$t_i - \tau \cdot j \leqslant T$ {\rm and} {\tt Flag}$_i = 0$}{
            Communicate $(t_i, t_i + \Delta_i)$ to client $i$\;
            Set $\texttt{Flag}_i = 1$\;
        }
    }

    \textbf{4.} Increase $j$ by 1\;
}
\caption{Dynamic-WOS procedure.}
\label{alg:online}
\end{algorithm}

\begin{remark}
\normalfont
Observe that in the regime that $\tau$ is small, the timing of the update and the associated time window can be precomputed. Specifically, when we would recompute the optimal (static) time windows for the remaining clients $N+1,\ldots,n$ upon completion of client $N$, we need the minimizer $(\boldsymbol{t}^{(N)}, \boldsymbol{\Delta}^{(N)}) \in \mathscr{U}_{n-N}$ of $\text{WOS}_{n - N}(\boldsymbol{B}[N], \omega, \mathscr{P})$, with ${\bs B}[N]:=(B_{N+1},\ldots,B_n)$.
Since $t_i^{(N)}$ decreases in $N$, it follows that if $T\in[t_i^{(N)} , t_i^{(N-1)})$, then the customer $i \ge N$ receives its update moment while the service provider travels from customer $N-1$ to $N$. 
It only remains to find the exact elapsed travel time of customer $N$ such that the start of the time window and the current time equals $T$. This moment can be determined in an iterative manner similar to Algorithm~\ref{alg:online}, for $\tau$ small.
\hfill $\Diamond$
\end{remark}

To show the achievable gains of a single timely update, we perform a set of numerical experiments in 
which we evaluate the performance of DWOS relative to our static scheduling procedure. Unless otherwise specified, an experiment consists of $10\,000$ simulation runs, and uses the following default parameters: $n = 25$ clients, a normally distributed service time with mean $\mu = 10$ and standard deviation $\sigma = 2.5$, advance notice threshold $T = 30$ minutes, and penalty parameters $\alpha = 0.1$, $\beta = 1.1$, and $\omega = 0.5$. After demonstrating
the impact of the parameters under such
a stylized setting, we also present the difference between WOS and DWOS in a
real-world context, using the last-mile delivery dataset \cite{Konovalenko2024LastMile}.

\begin{figure}
    \centering
\scalebox{0.5}{\begin{tikzpicture}
\begin{axis}[
    xmin=0,
    xlabel={Time (min)},
    ylabel={Client},
    ytick={1,...,25}, 
    xtick={0, 50, 100, 150, 200, 250},
    height=10cm,
    width=17cm,
    legend style={at={(0.97,0.18)}, anchor=north east, font=\small, fill=white},
    enlarge y limits=0.05,
    axis lines=box,
    tick align=outside,
    xtick pos=left,
    ytick pos=left,
]

\pgfplotstableread[col sep=comma]{Figures_latex/window_data_N1.csv}\datatable

\addlegendimage{line legend, thick, gray, line width=7pt}
\addlegendentry{Static}
\addlegendimage{line legend, thick, blue, line width=3pt}
\addlegendentry{DWOS}
\foreach \i in {1,...,24} {%
    \pgfplotstablegetelem{\i-1}{static_window_start}\of\datatable
    \edef\staticstart{\pgfplotsretval}
    \pgfplotstablegetelem{\i-1}{static_window_end}\of\datatable
    \edef\staticend{\pgfplotsretval}
    \addplot [thick, gray, line width=7pt, solid] coordinates {(\staticstart,\i) (\staticend,\i)};
}

\foreach \i in {1,...,24} {%
    \pgfplotstablegetelem{\i-1}{dwos_window_start}\of\datatable
    \edef\dwosstart{\pgfplotsretval}
    \pgfplotstablegetelem{\i-1}{dwos_window_end}\of\datatable
    \edef\dwosend{\pgfplotsretval}
    \addplot [thick, blue, line width=3pt, solid] coordinates {(\dwosstart,\i) (\dwosend,\i)};
}

\end{axis}
\end{tikzpicture}}
    \caption{The latest communicated time window for each client 
    under WOS (gray) and DWOS (blue),
    for a single simulation run
    in a homogeneous setting}.
    \label{fig:single_sim}
\end{figure}
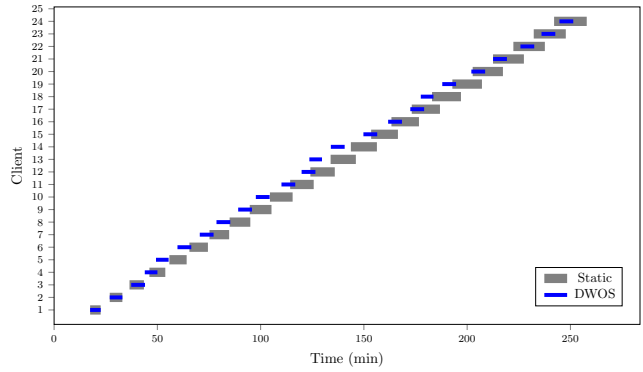

 First, to visually illustrate how the dynamically generated time windows under DWOS differ from the initial static schedule under WOS,
 Figure~\ref{fig:single_sim} 
 presents the latest communicated time windows to clients in a single simulation run.
  It can be observed that DWOS issues no updates for clients whose static windows begin before time~$T$ has elapsed, while sending a single update to clients whose static windows start after time $T$.
 In case the total realized travel times until that point are left (resp.\ right) from their mean, this will make the time windows start earlier (resp.\ later).
Moreover, as part of the travel times are no longer stochastic at the time of the update, there is less uncertainty, making the windows noticeably narrower under DWOS than under WOS. Also, given that all travel times in this example have the same distribution,
this window width becomes relatively uniform across clients. 

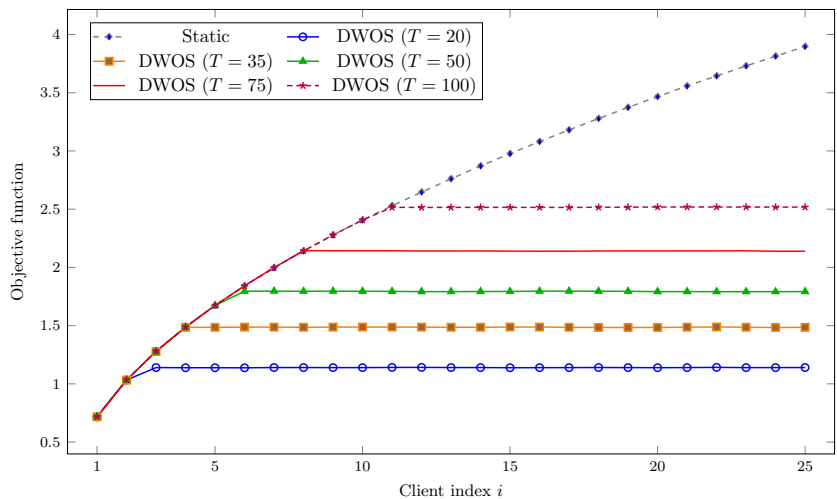
\begin{figure}
    \centering
    \scalebox{0.7}{\begin{tikzpicture}
\begin{axis}[
    width=16cm, height=10cm,
    xlabel={Client index $i$},
    ylabel={Objective function},
    legend style={ at={(0.03,0.97)},anchor=north west, scale=0.3, font=\small, fill=white, draw=black, legend columns=2, inner sep=1pt, outer sep=0pt, row sep=0pt, column sep=1ex,},
    xmin=0, xmax=26,
    xtick={1,5,10,15,20,25}
]
\addplot+[mark=diamond*, dashed, gray, thick]
    table[x=Client, y=Static_Total, col sep=comma] {Figures_latex/objective_comparison.csv};
\addlegendentry{Static}

\addplot+[mark=o, blue, thick]
    table[x=Client, y=DWOS_Total_T20, col sep=comma] {Figures_latex/objective_comparison.csv};
\addlegendentry{DWOS ($T=20$)}

\addplot+[mark=square*, orange, thick]
    table[x=Client, y=DWOS_Total_T35, col sep=comma] {Figures_latex/objective_comparison.csv};
\addlegendentry{DWOS ($T=35$)}

\addplot+[mark=triangle*, green!70!black, thick]
    table[x=Client, y=DWOS_Total_T50, col sep=comma] {Figures_latex/objective_comparison.csv};
\addlegendentry{DWOS ($T=50$)}

\addplot+[mark=triangle down, red, thick]
    table[x=Client, y=DWOS_Total_T75, col sep=comma] {Figures_latex/objective_comparison.csv};
\addlegendentry{DWOS ($T=75$)}

\addplot+[mark=star, purple, thick]
    table[x=Client, y=DWOS_Total_T100, col sep=comma] {Figures_latex/objective_comparison.csv};
\addlegendentry{DWOS ($T=100$)}
\end{axis}
\end{tikzpicture}}
    \caption{Client-wise total cost comparison between DWOS and static scheduling for various notice thresholds $T$.}
    \label{fig:stat_DWOS}
\end{figure}

\begin{figure}
    \centering
    \scalebox{0.7}{\begin{tikzpicture}
\begin{axis}[
    width=16cm, height=10cm,
    xlabel={Threshold $T$},
    ylabel={Relative Difference (\%)},
    legend style={
        at={(0.55,0.97)},
        anchor=north west,
        font=\small,
        fill=white,
        draw=black
    },
    grid=both,
    xtick={20,40,50,80,100,150,250},
    yticklabel style={/pgf/number format/.cd, fixed, precision=1},
    tick label style={font=\small},
    cycle list name=color list,
    every axis plot/.append style={thick, mark size=2.5pt},
    ymin=0, 
]

\addplot+[mark=o, thick, color=blue] table [
    x=T,
    y=alpha_0.1_omega_0.25_beta_1.3,
    col sep=comma
] {Figures_latex/rel_diff_total_plot.csv};
\addlegendentry{$\alpha=0.1$, $\omega=0.25$, $\beta=1.3$}

\addplot+[mark=square*, thick, color=orange] table [
    x=T,
    y=alpha_0.1_omega_0.5_beta_1.3,
    col sep=comma
] {Figures_latex/rel_diff_total_plot.csv};
\addlegendentry{$\alpha=0.1$, $\omega=0.5$, $\beta=1.3$}

\addplot+[mark=triangle*, thick, color=green!70!black] table [
    x=T,
    y=alpha_0.25_omega_0.5_beta_1.5,
    col sep=comma
] {Figures_latex/rel_diff_total_plot.csv};
\addlegendentry{$\alpha=0.25$, $\omega=0.5$, $\beta=1.5$}

\addplot+[mark=diamond*, thick, color=red] table [
    x=T,
    y=alpha_0.5_omega_0.25_beta_1.3,
    col sep=comma
] {Figures_latex/rel_diff_total_plot.csv};
\addlegendentry{$\alpha=0.5$, $\omega=0.25$, $\beta=1.3$}

\addplot+[mark=star, thick, color=purple] table [
    x=T,
    y=alpha_0.25_omega_0.75_beta_1.1,
    col sep=comma
] {Figures_latex/rel_diff_total_plot.csv};
\addlegendentry{$\alpha=0.25$, $\omega=0.75$, $\beta=1.1$}


\end{axis}
\end{tikzpicture}}
    \caption{Relative percentage cost difference between DWOS and static scheduling for various parameter combinations.}
    \label{fig:rel_diff}
\end{figure}
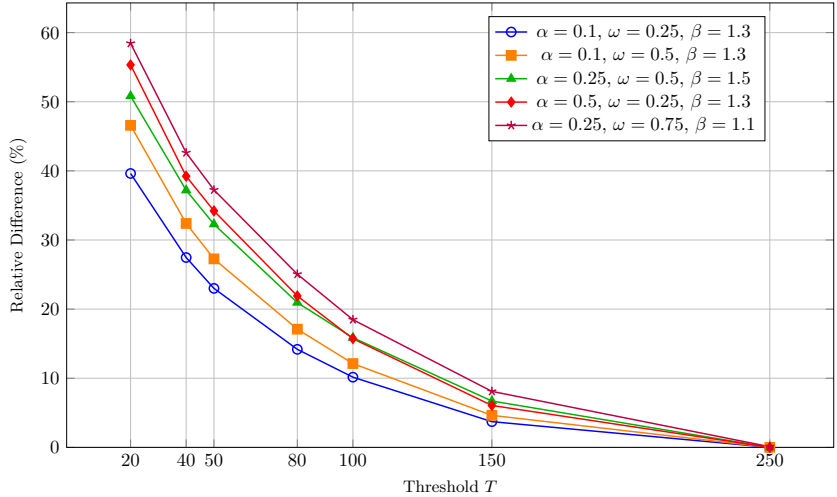

To illustrate the impact of updates on clients, Figure~\ref{fig:stat_DWOS} presents the average {\it per-client} costs across simulation runs for various thresholds $T$. The figure shows that DWOS consistently outperforms static scheduling by achieving lower client costs. For clients early in the schedule, the costs are identical across both methods, as DWOS does not issue updates to these clients either. However, while the costs under WOS increase significantly with the client’s position in the schedule, DWOS maintains a uniform cost profile across all clients.
Figure~\ref{fig:stat_DWOS} also shows that, naturally, the per-client cost difference between DWOS and WOS decreases as $T$ becomes smaller. This implies that a similar relation should hold for the total costs, which can be observed from Figure~\ref{fig:rel_diff},
which plots the relative difference of the total costs:
\begin{equation} \label{rel_gain}
\text{Relative Difference} = \frac{\text{Static Cost} - \text{DWOS Cost}}{\text{Static Cost}} \times 100\%.
\end{equation}
For all shown combinations of parameter values,
the costs decrease convexly in $T$. 
This pattern was consistently obtained across a wide range of parameter settings.

Table~\ref{tab:static_notice} further explores the implications of dynamic scheduling by quantifying the advance notice received by clients under DWOS. 
In this example, we take $T=30$ to determine when DWOS updates may be triggered. It is important to distinguish the DWOS threshold $T$ from the `notice thresholds' that clients experience; the latter are used solely to describe the distribution of the `advance notice', i.e., the time gap between the DWOS update and the start of the client's static time window. For selected client indices, we compute the proportion of simulation runs where the advance notice is less than $10$, $15$, or $25$ minutes. We observe that even for the 25th client, the advance notice is under 10 minutes in less than 10$\%$ of the cases.

\begin{table}[htbp]
\centering
{\small
\caption{Proportion of time the client receives static advance notice, with increasing client index. The values show the percentage of instances below time thresholds (in minutes), and the mean notice time.}
\label{tab:static_notice}
\begin{tabular}{lcccc}
\toprule
Client index & \% $<$ 10 min & \% $<$ 15 min & \% $<$ 25 min & Mean advance notice (min) \\
\midrule
5  & 0.00 & 0.00 & 5.97  & 28.90 \\
10 & 0.17 & 1.97 & 33.68 & 27.73 \\
15 & 2.22 & 8.46 & 42.73 & 26.32 \\
20 & 5.81 & 14.07 & 46.95 & 25.62 \\
25 & 9.46 & 18.91 & 49.73 & 24.34 \\
\bottomrule
\end{tabular}
}
\end{table}


To test how well DWOS performs beyond the studied
setting with normally distributed travel times that
are homogeneous over the clients, Figure~\ref{fig:boxDWOS}
shows the relative difference under heterogeneous settings. 
More specifically, the figure shows the result of $5\,000$ simulation runs, where in each run we randomly selected the system-wide parameters: similar to Figure~\ref{fig:boxplot_errors}, we randomly pick $\omega\in \{0.25, 0.5, 0.75\}$, $\alpha\in\{0.01, 0.05, 0.1, 0.15\}$, and $\beta\in\{1.1, 1.2, 1.3, 1.4, 1.5\}$. 
Moreover, for each client the travel-time distribution is randomly chosen among the normal, Lognormal, and Weibull families, with the mean $\mu\in\{10,12,15\}$ drawn at random.
As can be observed from Figure~\ref{fig:boxDWOS}, 
the use of DWOS generally significantly improves performance compared to static WOS, especially for smaller values of $T$. 
This confirms that DWOS also adapts well to heterogeneous client behavior, and that the impact of the updates is not just observed for homogeneous and normally distributed travel times.

\begin{figure}
    \centering
\scalebox{0.5}{\begin{tikzpicture}
\begin{axis}[
    xmin=0,
    xlabel={Time (min)},
    ylabel={Client},
    ytick={1,...,25}, 
    xtick={0, 200, 400, 600, 800, 1000},
    height=10cm,
    width=17cm,
    legend style={at={(0.97,0.18)}, anchor=north east, font=\small, fill=white},
    enlarge y limits=0.05,
    axis lines=box,
    tick align=outside,
    xtick pos=left,
    ytick pos=left,
]

\pgfplotstableread[col sep=comma]{Figures_latex/window_data_N_real.csv}\datatable

\addlegendimage{line legend, thick, gray, line width=7pt}
\addlegendentry{Static}
\addlegendimage{line legend, thick, blue, line width=3pt}
\addlegendentry{DWOS}
\foreach \i in {1,...,24} {%
    \pgfplotstablegetelem{\i-1}{static_window_start}\of\datatable
    \edef\staticstart{\pgfplotsretval}
    \pgfplotstablegetelem{\i-1}{static_window_end}\of\datatable
    \edef\staticend{\pgfplotsretval}
    \addplot [thick, gray, line width=7pt, solid] coordinates {(\staticstart,\i) (\staticend,\i)};
}

\foreach \i in {1,...,24} {%
    \pgfplotstablegetelem{\i-1}{dwos_window_start}\of\datatable
    \edef\dwosstart{\pgfplotsretval}
    \pgfplotstablegetelem{\i-1}{dwos_window_end}\of\datatable
    \edef\dwosend{\pgfplotsretval}
    \addplot [thick, blue, line width=3pt, solid] coordinates {(\dwosstart,\i) (\dwosend,\i)};
}

\end{axis}
\end{tikzpicture}}
    \caption{The latest communicated time window for each client 
    under WOS (gray) and DWOS (blue),
    for a single simulation run
    sampled from the data set \cite{Konovalenko2024LastMile}.}
    \label{fig:single_sim_2}
\end{figure}
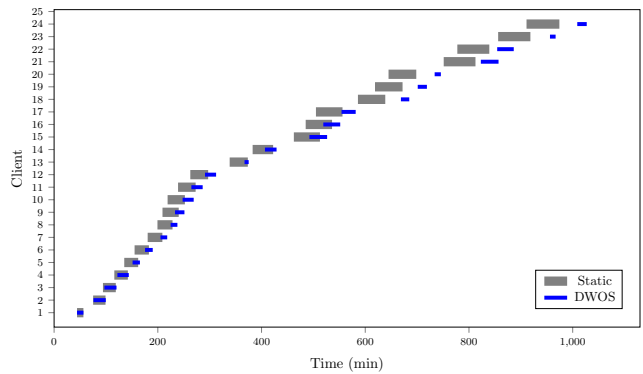

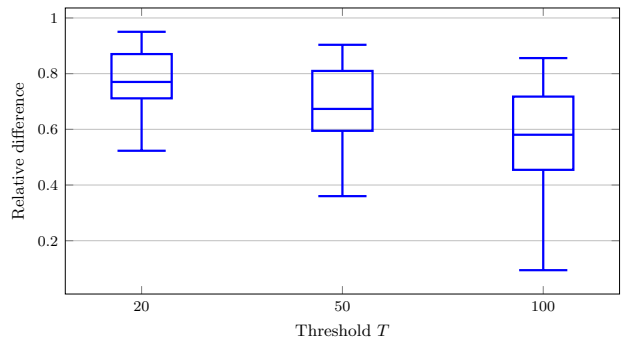
\begin{figure}
    \centering
    \scalebox{0.7}{\begin{tikzpicture}
\begin{axis}[
    width=12cm, height=7cm,
    ylabel={Relative difference},
    xlabel={Threshold $T$},
    xtick={1,2,3},
    xticklabels={20, 50, 100},
    ymajorgrids,
    boxplot/draw direction=y,
    boxplot/box extend=0.30,
    every box/.style={draw=none},
    /pgfplots/boxplot/every lower quartile/.append style={
        /pgfplots/draw/.append code={
            \draw[blue]
              (boxplot box cs:lower quartile,0) --
              (boxplot box cs:upper quartile,0);
        }
    },
]



\addplot+[
    boxplot prepared={
        median=0.77020679,
        upper quartile=0.87032764,
        lower quartile=0.7113482,
        upper whisker=0.95029091,
        lower whisker=0.52304875
    },
    draw=blue,
] coordinates {};



\addplot+[
    boxplot prepared={
        median=0.67335033,
        upper quartile=0.80980167,
        lower quartile=0.59468225,
        upper whisker=0.903633,
        lower whisker=0.3600798
    },
    draw=blue,
] coordinates {};



\addplot+[
    boxplot prepared={
        median=0.58057845,
        upper quartile=0.71759044,
        lower quartile=0.45454379,
        upper whisker=0.85590983,
        lower whisker=0.09415275
    },
    draw=blue,
] coordinates {};

\end{axis}
\end{tikzpicture}}
    \caption{Boxplot assessing the relative gain in DWOS for different thresholds in heterogeneous settings,
    in which the realization of $B_i$ is 
    a travel-time value from the data set 
    \cite{Konovalenko2024LastMile}, and
    the mean and standard deviation of $B_i$ 
    are
    fitted using this data set (see Appendix~\ref{app:data}). The whiskers correspond to the $5\%$- and $95\%$-percentiles. The outliers are omitted.}
    \label{fig:boxDWOS_figure_12}
\end{figure}

To further investigate the impact of heterogeneity in a real-world context, we repeated the above experiment using a last-mile delivery dataset \cite{Konovalenko2024LastMile}.
This dataset consists of trajectories of delivery vehicles and includes, among other variables, the time and Haversine distance between consecutive stops.
In short (for more detail, see
Appendix~\ref{app:data}), in each
simulation of the experiment, 
we sampled for
each client~$i$ a pair
of stops with corresponding (distance, time)-pair. 
The recorded time-value
is used as the realization
of the travel-time distribution
$B_i$, and a mixture of linear regression
models determines the mean and standard
deviation of $B_i$.
Figure~\ref{fig:single_sim_2}
displays the WOS and DWOS time windows
for a single simulation run for
25 sampled clients, again demonstrating that DWOS produces substantially narrower windows than WOS.
The resulting relative differences for
three different $T$ values and 1000 simulation runs
are presented in Figure~\ref{fig:boxDWOS_figure_12},
which shows that the conclusions from
Figure~\ref{fig:boxDWOS} extend to 
a practical setting in which the true
underlying travel time distributions are
unknown. It should be noted that the
distances in the data set are typically of
larger order than the selected means
in Figure~\ref{fig:boxDWOS} 
(see e.g.\ Figure~\ref{fig:single_sim_2}), such that
the relative difference between WOS
and DWOS is even present for the largest value of $T$.

Overall, these results demonstrate that DWOS consistently improves scheduling performance compared to the static WOS approach. The method proves robust across a variety of parameter settings and distributional assumptions, offering substantial gains even in heterogeneous and uncertain environments.


\begin{figure}
    \centering
    \scalebox{0.7}{\begin{tikzpicture}
\begin{axis}[
    width=12cm, height=7cm,
    ylabel={Relative difference},
    xlabel={Threshold $T$},
    xtick={1,2,3},
    xticklabels={20, 50, 100},
    ymajorgrids,
    boxplot/draw direction=y,
    boxplot/box extend=0.30,
    every box/.style={draw=none},
    /pgfplots/boxplot/every lower quartile/.append style={
        /pgfplots/draw/.append code={
            \draw[blue]
              (boxplot box cs:lower quartile,0) --
              (boxplot box cs:upper quartile,0);
        }
    },
]



\addplot+[
    boxplot prepared={
        median=0.55964732,
        upper quartile=0.70627131,
        lower quartile=0.39595682,
        upper whisker=0.80753565,
        lower whisker=0.07622513
    },
    draw=blue,
] coordinates {};



\addplot+[
    boxplot prepared={
        median=0.39565948,
        upper quartile=0.57574913,
        lower quartile=0.19482031,
        upper whisker=0.71127292,
        lower whisker=-0.14317722
    },
    draw=blue,
] coordinates {};



\addplot+[
    boxplot prepared={
        median=0.18531792,
        upper quartile=0.30045863,
        lower quartile=-0.06204421,
        upper whisker=0.49663288,
        lower whisker=-0.40069844
    },
    draw=blue,
] coordinates {};

\end{axis}
\end{tikzpicture}}
    \caption{Boxplot assessing the relative gain in DWOS for different thresholds in heterogeneous parameter settings,
    in which the underlying distribution of $B_i$ is either normal, Lognormal, or Weibull. The whiskers correspond to the $5\%$- and $95\%$-percentiles. The outliers are omitted.}
    \label{fig:boxDWOS}
\end{figure}
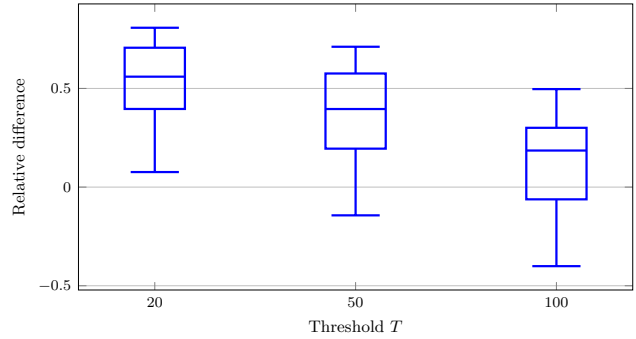

\section{Conclusive remarks} \label{sec: conclusive_remarks}

The objective of this paper was to develop a structured framework for designing and dynamically updating service time windows in delivery and appointment-based systems.
We considered a single-vehicle scenario with a given route and $n$ delivery points, where travel and service times are random. The first part of the paper addressed optimal time window assignment at the trip's start, the optimality criterion containing penalties for tardiness, lateness, and window width. We examined two variants: WOS (variable window widths) and UWOS (uniform window widths). For WOS,
it has been shown that determining optimal time windows reduces to one-dimensional root-finding problems, with explicit solutions for certain
specific distributions and an efficient approximation method for more general cases.
For UWOS, it is explained that solving an $n+1$-dimensional root-finding problem is required. Interestingly, although UWOS is often viewed as the fairer approach — since all customers receive equal window widths — numerical results revealed that the associated penalty costs are not evenly distributed among customers.

The second part of the paper studied the effect of
dynamically updating the initial time windows during the delivery trip.
To balance the need for timely communication with the risk of overwhelming customers with frequent updates, the proposed strategy sends a single notification once the remaining time until the start of the optimal window falls below a pre-specified threshold $T$.
Extensive numerical experiments demonstrated that these updated time windows are generally substantially narrower than the initially assigned ones. Furthermore, when the update threshold is small, the improvement in the objective function is significant.

There are several directions for follow-up research
that emerge from this work. 
For example, one could consider a large-scale empirical validation of the presented framework using real-world data.
Another direction involves designing a framework in which requests appear dynamically and randomly before the delivery trip, but must be assigned a time window on the spot (or before some cut-off point, as in \cite{CSMW}). While this paper provides a numerical comparison between  UWOS and WOS, it would also be interesting to establish theoretical results for their difference.
Finally, a more advanced dynamic updating procedure could be explored, in which customers receive multiple updates or are notified only when there is a significant change between the current and previously communicated time window.


\bibliographystyle{abbrv}
\bibliography{references}

\begin{thebibliography}{10}

\bibitem{ahmadi2017}
A.~Ahmadi-Javid, Z.~Jalali, and K.~J. Klassen.
\newblock Outpatient appointment systems in healthcare: A review of optimization studies.
\newblock {\em European Journal of Operational Research}, 258(1):3--34, 2017.

\bibitem{bekker2023}
R.~Bekker, B.~Bharti, L.~Lan, and M.~Mandjes.
\newblock A queueing-based approach for integrated routing and appointment scheduling.
\newblock {\em European Journal of Operational Research}, 318(2):534--548, 2024.

\bibitem{BKCE}
M.~Burian, C.~K\"ohler, M.~Campbell, and J.~Ehmke.
\newblock Service time window selection for attended home deliveries: a case study for urban and rural areas.
\newblock {\em Central European Journal of Operations Research}, 32:267--294, 2024.

\bibitem{CSMW}
S.~Celik, A.~Schrotenboer, L.~Martin, and T.~van Woensel.
\newblock {Is Waiting Worth It? The Value of Delaying Time Window Assignment in Vehicle Routing Problems}.
\newblock {\em Transportation Research Part B: Methodological}, 204:103381, 2026.

\bibitem{dalmeijer2018}
K.~Dalmeijer and R.~Spliet.
\newblock A branch-and-cut algorithm for the time window assignment vehicle routing problem.
\newblock {\em Computers \& Operations Research}, 89:140--152, 2018.

\bibitem{DSW}
K.~Dalmeijer, R.~Spliet, and A.~P. Wagelmans.
\newblock Dynamic time window adjustment.
\newblock Technical Report EI2019-22, Econometric Institute, Erasmus School of Economics, Erasmus University of Rotterdam, 2019.

\bibitem{de2021}
M.~{de Kemp}, M.~Mandjes, and N.~Olver.
\newblock Performance of the smallest-variance-first rule in appointment sequencing.
\newblock {\em Operations Research}, 69(6):1909--1935, 2021.

\bibitem{hassin2008scheduling}
R.~Hassin and S.~Mendel.
\newblock Scheduling arrivals to queues: A single-server model with no-shows.
\newblock {\em Management Science}, 54(3):565--572, 2008.

\bibitem{HADJ}
M.~Hoogeboom, Y.~Adulyasak, W.~Dullaert, and J.~Patrick.
\newblock The robust vehicle routing problem with time window assignments.
\newblock {\em Transportation Science}, 55(2):395--413, 2021.

\bibitem{hosseini2025}
D.~Hosseini, B.~Rostami, and M.~Araghi.
\newblock Service time window design in last-mile delivery.
\newblock {\em ArXiv}, page \url{https://arxiv.org/abs/2508.01032}, 2025.

\bibitem{JAB}
O.~Jabali, R.~Leus, T.~{Van Woensel}, and T.~{de Kok}.
\newblock Self-imposed time windows in vehicle routing problems.
\newblock {\em OR Spectrum}, 37(2):331--352, 2013.

\bibitem{kemper2014}
B.~Kemper, C.~Klaassen, and M.~Mandjes.
\newblock Optimized appointment scheduling.
\newblock {\em European Journal of Operational Research}, 239(1):243--255, 2014.

\bibitem{KEC}
C.~K\"ohler, J.~Ehmke, and A.~Campbell.
\newblock Flexible time window management for attended home deliveries.
\newblock {\em Omega}, 91:102023, 2020.

\bibitem{Konovalenko2024LastMile}
A.~Konovalenko, L.~M. Hvattum, and K.~A.~H. Iversen.
\newblock Last-mile delivery route deviations dataset: Planned vs. actual routes.
\newblock \url{https://doi.org/10.17632/kkwgfvmtxn.1}, 2024.
\newblock Mendeley Data.

\bibitem{KONOVALENKO2026129921}
A.~Konovalenko, L.~M. Hvattum, and K.~A.~H. Iversen.
\newblock Predicting last-mile delivery route deviations using machine learning.
\newblock {\em Expert Systems with Applications}, 298:129921, 2026.

\bibitem{LLM}
G.~Laporte, F.~Louveaux, and H.~Mercure.
\newblock The vehicle routing problem with stochastic travel times.
\newblock {\em Transportation Science}, 26(3):161--170, 1992.

\bibitem{NPGAA}
F.~Neves-Moreira, D.~Pereira~da Silva, L.~Guimar\~aes, P.~Amorim, and B.~Almada-Lobo.
\newblock The time window assignment vehicle routing problem with product dependent deliveries.
\newblock {\em Transportation Research Part E: Logistics and Transportation Review}, 116:163--183, 2018.

\bibitem{PRLD}
R.~Paradiso, R.~Roberti, D.~Lagan\'a, and W.~Dullaert.
\newblock An exact solution framework for multitrip vehicle-routing problems with time windows.
\newblock {\em Operations Research}, 68(1):180--198, 2020.

\bibitem{PR}
C.~D. Pegden and M.~Rosenshine.
\newblock Scheduling arrivals to queues.
\newblock {\em Computers \& Operations Research}, 17(4):343--348, 1990.

\bibitem{spliet2015}
R.~Spliet and A.~F. Gabor.
\newblock The time window assignment vehicle routing problem.
\newblock {\em Transportation Science}, 49(4):721--731, 2015.

\bibitem{UGM}
M.~Ulmer, J.~Goodson, and B.~Thomas.
\newblock Optimal service time windows.
\newblock {\em Transportation Science}, 58(2):394--411, 2024.

\bibitem{VRT}
A.~Vareias, P.~Repoussis, and C.~Tarantilis.
\newblock Assessing customer service reliability in route planning with self-imposed time windows and stochastic travel times.
\newblock {\em Transportation Science}, 53(1):256--281, 2017.

\bibitem{wang1997}
P.~Wang.
\newblock Optimally scheduling {N} customer arrival times for a single-server system.
\newblock {\em Computers \& Operations Research}, 24(8):703--716, 1997.

\bibitem{WAS}
K.~Wa{\ss}muth, C.~K\"ohler, N.~Agatz, and M.~Fleischmann.
\newblock Demand management for attended home delivery—a literature review.
\newblock {\em European Journal of Operational Research}, 311(3):801--815, 2023.

\bibitem{YSBS}
X.~Yu, S.~Shen, B.~Badri-Koohi, and H.~Seada.
\newblock Time window optimization for attended home service delivery under multiple sources of uncertainties.
\newblock {\em Computers \& Operations Research}, 150:106045, 2023.

\bibitem{zhan2021}
Y.~Zhan, Z.~Wang, and G.~Wan.
\newblock Home service routing and appointment scheduling with stochastic service times.
\newblock {\em European Journal of Operational Research}, 288(1):98--110, 2021.

\end{thebibliography}
\appendix

\section{Proofs of Section~\ref{sec:wos}}
\label{app: proofs_wos}

{\it Proof of Theorem~\ref{THM1}.}
 Realize that $g(\Delta):={\mathscr P}'(\Delta)$ is increasing in $\Delta$, and assume that $g(\cdot)$ is invertible. Hence,
\[\Delta_i(t_i) = g^{-1}((1-\omega)F_i(t_i)),\]
which is increasing in $t_i$ (as it is the composition of the increasing functions $g^{-1}(\cdot)$ and $F_i(\cdot)$).
We thus obtain the equation
\[{\mathfrak f}^-_i(t_i):= \omega\left(1-F_i\left(t_i +g^{-1}((1-\omega)F_i(t_i)\right)\right) = (1-\omega)F_i(t_i) =:{\mathfrak f}^+_i(t_i).  \]
Observe that ${\mathfrak f}^-_i(\cdot)$ is decreasing, with ${\mathfrak f}^-_i(0)= \omega$ and ${\mathfrak f}^-_i(\infty)=0$, and that ${\mathfrak f}^+_i(\cdot)$ is increasing, with ${\mathfrak f}^+_i(0)= 0$ and ${\mathfrak f}^+_i(\infty)=1-\omega.$ Conclude that there is a unique positive root for any $\omega\in(0,1).$ 

To verify whether the found stationary point is indeed a {\it minimum}, we consider the Hessian, using the familiar notation $f_i(t):=F_i'(t)$,
\[ \left(\begin{array}{cc}(1-\omega)f_i(t_i)+\omega f_i(t_i+\Delta_i)&\omega f_i(t_i+\Delta_i)\\
\omega f_i(t_i+\Delta_i)
&\omega f_i(t_i+\Delta_i)+{\mathscr P}''(\Delta_i)\end{array}\right).\]
This matrix is diagonally-dominant, hence positive definite (where it is used that ${\mathscr P}(\cdot)$ being strictly convex entails ${\mathscr P}''(\Delta_i)>0$). It follows that it has positive eigenvalues, which implies that any stationary point (i.e., solution of the first-order conditions) is a minimum.
\hfill$\Box$

{\it Proof of Lemma~\ref{lem: symmetry}.}
By Theorem~\ref{THM1}, $(t(\omega'), \Delta)$ is the unique
solution to \eqref{FOC1} and \eqref{FOC2}
for $\omega = \omega'$, such that
\begin{align} \label{eq: opt_omega_prime}
    (1-\omega') F_i(t(\omega')) = \omega'\left(1-F_i(t(\omega') + \Delta )\right) = \mathscr{P}'(\Delta).
\end{align}
Now, the symmetry of the probability density function of $S_i$ in $\tilde \mu_i$
yields that $F_i(\tilde \mu_i - x) = 1-F_i(\tilde \mu_i + x)$, 
from which we obtain
\begin{align*}
    F_i(t(1-\omega'))
    =F_i\left(\tilde\mu_i - (\tilde\mu_i-(1-\omega'))\right)
    =1-F_i(-t(1-\omega') + 2\tilde\mu_i)
    =1-F_i(t(\omega') + \Delta)
\end{align*}
and
\begin{align*}
    1-F_i(t(1\!-\!\omega') + \Delta)
    =1-F_i(\tilde\mu_i + t(1\!-\!\omega') + \Delta - \tilde\mu_i)
    =F_i(2\tilde\mu_i - t(1\!-\!\omega') - \Delta)
    =F_i(t(\omega')).
\end{align*}
Hence, \eqref{eq: opt_omega_prime} yields that
both $\omega' F_i(t(1-\omega'))$ and
$(1-\omega')\left(1-F_i(t(1-\omega') + \Delta)\right)$
equal $\mathscr{P}'(\Delta)$, making
$[t(1-\omega'), t(1-\omega') + \Delta]$ the optimal
window for $\omega = 1-\omega'$. Note that
for $\omega = \frac{1}{2}$, \eqref{eq: omega_prime} implies
\begin{align*}
    \Bigg|t\left(\frac{1}{2}\right) + \Delta_i - \tilde\mu_i\Bigg|
    = \Bigg|t\left(\frac{1}{2}\right) - \tilde\mu_i\Bigg|,
\end{align*}
as desired.
\hfill$\Box$

\section{Mean and variance of residual travel times}\label{app:cond}
It is observed that in the dynamic schedule {\it residual} travel times play a crucial role; see specifically Eqn.~\eqref{eq:conditional}. In the spirit of what we discussed in the paper, we consider in this appendix the case in which we replace the residual travel times by their normally distributed counterparts. This leaves us with the task of evaluating their means and variances.

Recall that the residual travel time, conditional on it being larger than $b$, is given via 
\[{\mathbb P}(B_{j}^{({\rm res})}>t) := {\mathbb P}(B_{j}>b+t\,|\,B_j>b)=\frac{{\mathbb P}(B_{j}>t+b)}{{\mathbb P(B_{j}>b)}}.\]
Now suppose that $B_j$ has a normal distribution with mean $\mu_j$ and variance $\sigma_j^2.$ Define $Z:=(B_j-\mu_j)/\sigma_j$, with $Z$ standard normal. 
Then, 
\[{\mathbb E}[B_j\,|\, B_j>b] = {\mathbb E}[\mu_j+\sigma_j Z \,|\, B_j>b] = \mu_j + \sigma_j\,{\mathbb E}[Z\,|\, B_j>b] = \mu_j + \sigma_j\, {\mathbb E}[Z \,|\, Z > z]\]
where $z$ is the centered and normalized version of $b$, i.e., $z\equiv z_b:=(b-\mu_j)/\sigma_j$.
Note that, with $\varphi(\cdot)$ the standard normal probability density function and $\Phi(z)$ the corresponding cumulative distribution function, for any $z$,
\[{\mathbb E}[Z \,|\, Z >  z] = \frac{1}{1 - \Phi(z)} \int_{z}^{\infty} v \varphi(v) \, {\rm d}v.\]
Computing the integral, we readily find \[\bar{\mu}_j(b) :={\mathbb E}[B_j^{\rm (res)}]={\mathbb E}[B_j\,|\, B_j>b]= \mu_j + \sigma_j \frac{\varphi(z_b)}{1 - \Phi(z_b)}.\]
We proceed by computing the conditional variance. The starting point is the definition
\[{\mathbb V}{\rm ar}\,(B_j \,|\, B_j> b) = {\mathbb E}[B_j^2 \,|\, B_j > b] - \bigl( {\mathbb E}[B_j \,|\, B_j > b] \bigr)^2.\]
Following the same reasoning as above, 
\[{\mathbb E}[B_j^2 \,|\, B_j > b] 
= \mu_j^2 + 2\mu_j\sigma_j\, {\mathbb E}[Z \,|\, Z > z] + \sigma_j^2 \,{\mathbb E}[Z^2 \,|\, Z > z],\]
where
\[{\mathbb E}[Z^2 \,|\, Z > z] = \frac{1}{1 - \Phi(z)}\int_{z}^{\infty} v^2 \varphi(v) \, {\rm d}v=\frac{z \varphi(z) + 1 - \Phi(z)}{1 - \Phi(z)}=1 + \frac{z \varphi(z)}{1 - \Phi(z)}.\]
Noting that 
\begin{align*}{\mathbb E}[B_j^2 \,|\, B_j > b] &= \mu_j^2 + 2\mu_j\sigma_j \frac{\varphi(z)}{1 - \Phi(z)} + \sigma_j^2 \left( 1 + \frac{z \varphi(z)}{1 - \Phi(z)} \right),\\
 \bigl( {\mathbb E}[B_j \,|\, B_j > b] \bigr)^2 &= \mu_j^2 + 2\mu_j\sigma_j \frac{\varphi(z)}{1 - \Phi(z)} + \sigma_j^2 \left( \frac{\varphi(z)}{1 - \Phi(z)} \right)^2,\end{align*}
we thus obtain, after elementary computations,
\[\bar{\sigma}_j^2(b) := {\mathbb V}{\rm ar}\,(B_j^{\rm (res)})={\mathbb V}{\rm ar}\,(B_j \,|\, B_j> b)=\sigma_j^2 \left( 1 + \frac{z_b \varphi(z_b)}{1 - \Phi(z_b)} - \left( \frac{\varphi(z_b)}{1 - \Phi(z_b)} \right)^2 \right).\]


\section{Travel time data}\label{app:data}
In Figure~\ref{fig:boxDWOS_figure_12}, we have applied the proposed online methodology to a real-world setting using last-mile delivery data from a logistics company \cite{Konovalenko2024LastMile}. In this appendix we provide more background on the underlying experiments. The dataset, originally introduced in \cite{KONOVALENKO2026129921}, contains trajectories of delivery vehicles. For each pair of consecutive stops along a trajectory, the dataset provides both the recorded elapsed time and the corresponding Haversine distance between the stops. The recorded time captures not only travel time but also service-related activities, such as unloading the vehicle and handing over parcels.

For this study, we collected all 
(distance, time) pairs with a nonzero recorded time and a distance of at least 2~kilometers. To mitigate the effects of noise and potential data errors, we cleaned the dataset by excluding pairs for which the ratio of time to distance fell within the top or bottom 2.5$\%$ of the distribution. From the remaining 81911 entries, 70$\%$ were used for training, while the remaining 30$\%$ formed the sample set used in the experiment shown in Figure~\ref{fig:boxDWOS_figure_12}. 
The latter will be referred to
as `test set'.

\begin{figure}
    \centering
    \includegraphics[width=0.5\linewidth]{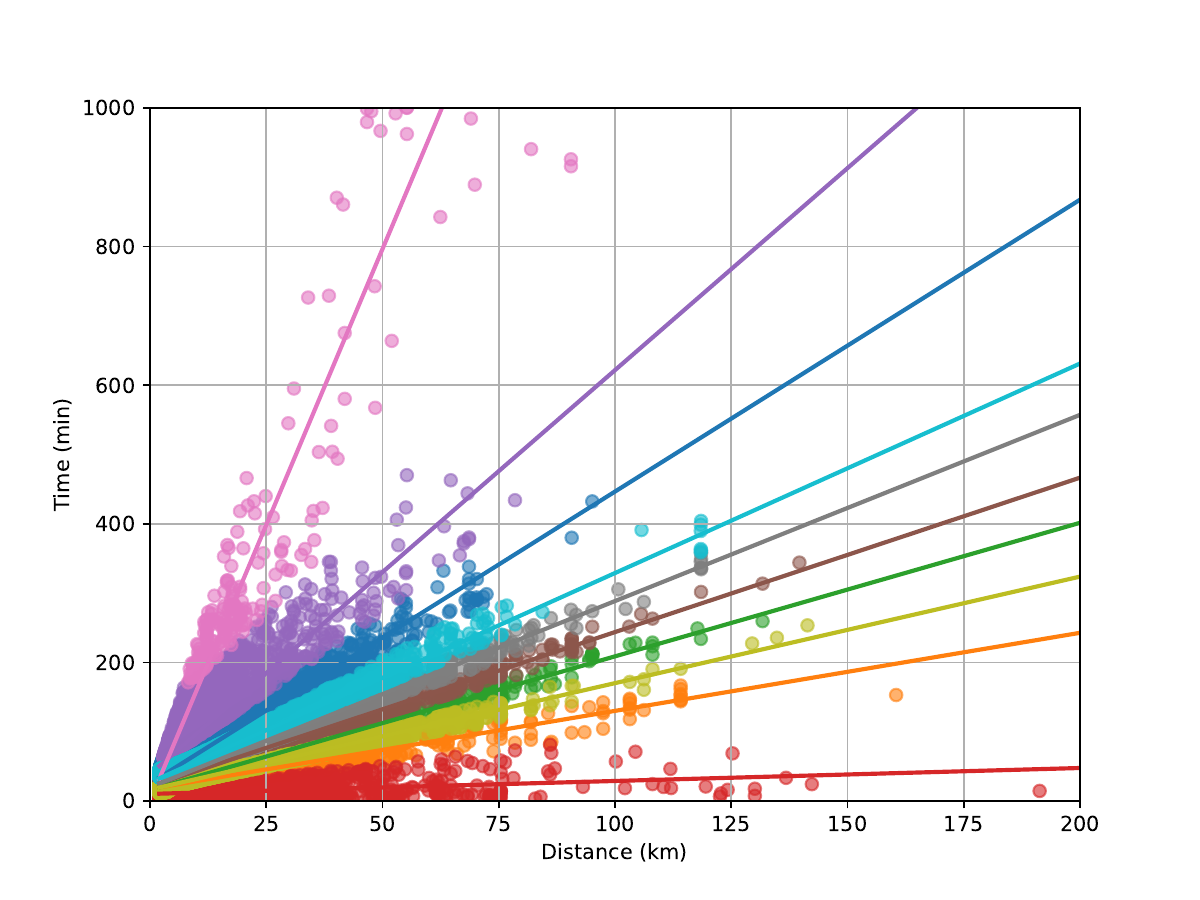}
    \caption{Mixture regression model}
    \label{fig:scatter_plot_data}
\end{figure}

Given the anticipated
linear relation between distance and travel time,
a natural candidate model is a 
linear regression model with Gaussian noise.
However, the data shows that a single
linear line would not 
fit the (distance, time) pairs
well.
This is, for example, explained by the fact that there
is a substantial amount of heterogeneity between different GPS
pairs. Indeed, depending on the
network topology, the Haversine distance may not always
represent the true traveled distance well,
and the traffic conditions and speed limits
vary over the 
different roads in the road network as well. 
To account for these differences,
we use a \textit{mixture} of ten linear regression models with Gaussian noise, in which
the $k$-th model has intercept
$a_k \geq 0$, slope 
$b_k \geq 0$, and
standard deviation $\sigma_k$:
\begin{align*}
    \text{ time } =
    a_k + b_k \cdot \text{ distance }
    + \epsilon, 
    \quad \quad 
    \epsilon \sim \mathcal{N}(0, \sigma_k^2).
\end{align*}
The non-negativity of $a_k$ and $b_k$ reflects the non-negativity of time and distance as model variables. To guarantee this property in the estimation procedure, the parameters $a_k$, $b_k$, and $\sigma_k$ are estimated 
{with the training data
via the EM algorithm},
where inequality constraints are enforced in the M-step to ensure non-negativity of the resulting coefficients $a_k$ and $b_k$.
{Moreover,
to keep the noise small,
$\sigma_k$ is enforced
to be at most $3\bar{y}$, 
where $\bar{y}$ denotes
the mean distance of the training data}.
To initialize the EM-algorithm, $a_k$
is chosen uniformly between $0$ and 
$\bar{y}$,
$b_k$ is chosen uniformly between 0 and 
{the
absolute value of the mean of $(\text{time}_i - \bar{y})/\text{distance}_i$, $\sigma_k$ is set
as the standard deviation of the training
data distance}, and the weights of
each of the 10 models are
set 1/10. The resulting mixture model
is shown in Figure~\ref{fig:scatter_plot_data}.

Now, for each (distance, time) pair in the test set, a {\it maximum a posteriori estimator} is used to assign the observation to one of the ten linear models. This concretely means that, given the mixture weights of the models and the conditional Gaussian distribution of the time variable given the distance within each model, we compute the posterior probability that the observation was generated by each model; the observation is assigned to the model with the largest posterior probability.
Then,
if a (distance, time) pair from the test set is sampled to
represent client~$i$, 
the recorded value time
is used as the realization
of the travel time distribution
$B_i$.
If the pair is assigned to model~$k$, the value $a_k + b_k \cdot \text{distance}$ is used as the mean of $B_i$, and the standard deviation of $B_i$ is set to $\sigma_k$.

 \end{document}